\documentclass[onefignum,onetabnum,review]{siamart171218}
\usepackage{graphicx}
\usepackage{diagbox}
\usepackage{cite}
\usepackage{hyperref}
\usepackage{tablefootnote}

\usepackage{amsmath,amssymb,amsfonts}
\usepackage{algpseudocode}
\usepackage{algorithmicx}
\usepackage{multirow}
\usepackage{lineno}
\nolinenumbers
\newcommand{\eps}{\varepsilon}
\newcommand{\tv}{\bar{V}}
\newcommand{\tq}{\bar{Q}}

\newcommand{\vect}{\mathbf}
\newcommand{\ignore}[1]{}

\title{Iterated Gauss-Seidel GMRES }

\author{
Stephen~Thomas\thanks{National Renewable Energy Laboratory, Golden, CO}
\and Erin Carson\thanks{Charles University, Prague, CZ}
\and Miroslav Rozlo{\v{z}}n{\'{\i}}k\thanks{Czech Academy of Sciences, Institute of Mathematics, Prague, CZ.  Supported by Czech Academy of Sciences  (RVO 67985840) and by the Grant Agency
of the Czech Republic, Grant No. 23-06159S.} 
 \and Arielle Carr\thanks{Lehigh University, Bethlehem, PA}
\and Kasia \'{S}wirydowicz\thanks{Pacific Northwest National Laboratory, Richland, WA}
}

\begin{document}

\maketitle

\begin{abstract}
The GMRES algorithm of Saad and Schultz (1986) is an 
iterative method for approximately solving linear systems $A\vect{x}=\vect{b}$,
with initial guess $\vect{x}_0$ and residual $\vect{r}_0 = \vect{b} - A\vect{x}_0$. 
The algorithm employs the Arnoldi process to generate the Krylov basis vectors 
(the columns of $V_k$). It is well known that this process can be viewed as a  
$QR$ factorization of the matrix $B_k = [\: \vect{r}_0, AV_k\:]$  at each iteration. 
Despite an ${\cal O}(\eps)\kappa(B_k)$ loss of orthogonality, for unit roundoff $\eps$ and condition number $\kappa$, 
the modified Gram-Schmidt 
formulation was shown to be backward stable in the seminal paper by Paige
et al.~(2006).  
We present an iterated Gauss-Seidel formulation of
the GMRES algorithm (IGS-GMRES) based on the ideas of Ruhe (1983) and \'{S}wirydowicz et al. (2020).
IGS-GMRES maintains orthogonality to the level ${\cal O}(\eps)\kappa(B_k)$
or ${\cal O}(\eps)$,  depending on the choice of one or two iterations; 
for two Gauss-Seidel iterations, the computed Krylov basis vectors remain
orthogonal to working precision and  
the smallest singular value of $V_k$ remains close to one. 
The resulting GMRES method is thus backward stable. 
We show that IGS-GMRES can be implemented with only a single synchronization 
point per iteration, making it relevant to large-scale parallel computing environments. 
We also demonstrate that, unlike MGS-GMRES, in IGS-GMRES the relative Arnoldi residual
corresponding to the computed approximate solution
no longer stagnates above machine precision even for highly non-normal systems. 

\end{abstract}
   
\section{Introduction}

We consider linear systems of the form $A\vect{x}=\vect{b}$, where $A$ is an 
$n \times n$ real-valued matrix, solved via the generalized minimal residual method (GMRES) \cite{Saad86}. 
Let $\vect{r}_0 = \vect{b} - A\vect{x}_0$ 
denote the initial residual for the approximate solution $\vect{x}_0$.
Inside GMRES, the Arnoldi-$QR$ algorithm is applied to orthogonalize the
basis vectors for the Krylov subspace ${\cal K}_k(A,\vect{r}_0)$ 
spanned by the columns of the $n\times k$ matrix $V_k$, where $k \ll n$.
After $k$ iterations, this produces the $(k+1)\times k$ upper Hessenberg matrix 
$H_{k+1,k}$ in the Arnoldi expansion such that%
\begin{equation}
v_1= \vect{r}_0 / \rho, \quad  
     \rho := \|\vect{r}_0\|_2, \quad AV_k = V_{k+1}H_{k+1,k}.
\label{eq:arnoldi}
\end{equation}
The expansion (\ref{eq:arnoldi}) is equivalent to a 
$QR$ factorization of  $B_k = [\:\vect{r}_0,\:AV_k\:] = V_{k+1}
\left[ \: \rho\:\vect{e}_1 , \: H_{k+1,k} \: \right]$ and the columns of $V_{k}$
form an orthonormal basis for the Krylov subspace ${\cal K}_{k}(\:A,\:\vect{r}_0\:)$ 
\cite{2006--simax--paige-rozloznik-strakos}. The approximate solution in step $k$ is given
by $\vect{x}_k = \vect{x}_0 + V_k\vect{y}_k$,
where the vector $\vect{y}_k $ minimizes the Arnoldi residual
\begin{equation}
\| \:\rho\: \vect{e}_1 - H_{k+1,k}\vect{y}_k\:\|_2  
= \min_{y} \:
\| \: \rho\: \vect{e}_1 - H_{k+1,k}\vect{y}\:\|_2.
\label{eq:llse}
\end{equation}

When the Krylov vectors are orthogonalized via the 
modified Gram-Schmidt (MGS) algorithm in  finite-precision arithmetic, 
their loss of orthogonality is related in a 
straightforward way to the convergence of GMRES. Orthogonality 
among the Krylov vectors is effectively maintained until the norm-wise relative
backward error approaches the machine precision as discussed
in the work of Paige and Strako\v{s} \cite{Paige2002} and
Paige et al.~\cite{2006--simax--paige-rozloznik-strakos}.
The growth of the condition number $\kappa(B_k) = \sigma_{\max} (B_k)/
\sigma_{\min}(B_k)$, where $\sigma_{\max}(B_k)$ and $\sigma_{\min}(B_k)$ 
are the maximum and minimum singular values of the matrix $B_k$, respectively, 
is related to the 
 norm-wise backward error
\begin{equation}
\beta(\vect{x}_k) := \frac{\|\vect{b} - A\vect{x}_k\|_2}{\|\vect{b}\|_2 + \|A\|_2\|\vect{x}_k\|_2},
\quad k=1, 2, \ldots
\label{eq:nrbe}
\end{equation}
%
and it is observed that, in exact arithmetic, 
$\beta(\vect{x}_k) \: \kappa( B_k ) = 
{\cal O}(1)$.
Also in exact arithmetic, the orthogonality of the columns
implies the linear independence of the Krylov basis vectors. 

However, in 
finite-precision arithmetic, the columns of the computed $\tv_k$ are 
no longer orthogonal, as measured by $\|I - \tv_k^T \tv_k\|_F$, and
may deviate substantially from the level of machine precision, ${\cal O}(\eps)$.  
When linear independence is completely lost, the Arnoldi relative residual,
corresponding to the computed approximate solution $\vect{\bar{x}}_k$, 
stagnates at a certain level above ${\cal O}(\eps)$. 
This occurs when $\|S_k\|_2 = 1$, where 
\[
S_k = (I+L_k^T)^{-1}L_k^T,
\] 
and 
$L_k$ is the $k \times k$ strictly lower triangular part of $\tv_k^T\tv_k = I + L_k + L_k^T$. 
Note also that $I - S_k=(I+L_k^T)^{-1}$.
However, Paige et al. \cite{2006--simax--paige-rozloznik-strakos} demonstrated
that despite ${\cal O}(\eps) \kappa(B_k)$ loss of orthogonality, 
MGS-GMRES is backward stable for the solution of linear systems. 

Both modified and classical Gram-Schmidt versions with 
delayed re-orthogonalization (DCGS-2) were derived in 
\'{S}wirydowicz et al.~\cite{2020-swirydowicz-nlawa}, 
and these have been empirically observed to result in a backward stable GMRES method. 
The development of low-synchronization Gram-Schmidt and generalized 
minimal residual algorithms by  \'{S}wirydowicz et al.~\cite{2020-swirydowicz-nlawa}  
and  Bielich et al.~\cite{DCGS2Arnoldi} was largely 
driven by applications that need stable, yet scalable, solvers. 
An inverse compact $WY$ MGS algorithm is 
presented in~\cite{2020-swirydowicz-nlawa} and is 
based upon the application of the approximate projector
\[
P_k^{(1)} = I - \tv_k \: T_k^{(1)} \: \tv_k^T, \quad 
T_k^{(1)} \approx ( \tv_k^T\tv_k )^{-1},
\]
where $\tv_k$ is again $n\times k$, $I$ is the identity matrix of 
dimension $n$,
and $T_k^{(1)}$ is a $k\times k$ 
lower triangular correction matrix.
To obtain a low-synchronization algorithm requiring only a single global 
reduction per iteration, the normalization
is delayed to the next iteration. 
The matrix $T_k^{(1)} = (\:I + L_k\:)^{-1}$ is obtained from 
the strictly lower triangular part of $\tv_k^T\tv_k$, again denoted $L_k$.  Note that 
because $\tv_k$ has almost orthonormal columns, the norm of $L_k$ is small, 
and $T_k^{(1)}$ is nearly the identity matrix of dimension $k$. 
The last row of $L_{k}$ is constructed from the matrix-vector product
$L_{k,1:k-1} = \bar{\vect{v}}_{k}^T\tv_{k-1}$.

The purpose of the present work is to derive an 
iterated Gauss-Seidel formulation of the GMRES algorithm
based on the approximate solution of the normal
equations in the Gram-Schmidt projector,  as
described by Ruhe \cite{Ruhe83}, and the low-synchronization 
algorithms introduced by \'{S}wirydowicz et al. \cite{2020-swirydowicz-nlawa}.
$T_k^{(1)}$ represents the inverse compact $WY$ form of the MGS projector $P_k^{(1)}$
and corresponds to one Gauss-Seidel iteration for the approximate
solution of the normal equations
\begin{equation}\label{eq:normaleqn}
\tv_{k-1}^T \tv_{k-1} \vect{r}_{1:k-1,k} = \tv_{k-1}^T\:A\:\bar{\vect{v}}_{k}.
\end{equation}

Ruhe \cite[pg. 597]{Ruhe83} suggested applying
LSQR, whereas Bj\"orck \cite[pg. 312]{Bjorck94} recommended conjugate gradients
to solve the normal equations.
Instead, we apply two Gauss-Seidel iterations 
and demonstrate that the loss of orthogonality is at 
the level of ${\cal O}(\eps)$ 
{\it without any need for explicit reorthogonalization} in the Arnoldi-$QR$ 
algorithm. Resembling the Householder transformation-based 
implementation, the GMRES method with two Gauss-Seidel iterations 
is backward stable, where the 
stability result of Drko\v{s}ov\'{a} et al.~\cite{drsb1995} applies.
Even for extremely ill-conditioned and non-normal matrices, 
the two iteration Gauss-Seidel GMRES formulation matches
the Householder (HH)-GMRES of Walker \cite{Walker88},
where the Arnoldi residual of the
least square problem (\ref{eq:llse}) continues to
decrease monotonically without any stagnation.

\subsection*{Contributions} In this paper, we present an
iterated Gauss Seidel formulation of the GMRES algorithm \cite{Saad86}, and establish the backward 
stability of this approach. 
The computed Krylov basis vectors maintain orthogonality 
to the level of machine precision 
by projection onto their orthogonal complement and applying the 
low-synchronization Gram-Schmidt 
algorithms of  \'{S}wirydowicz et al.~\cite{2020-swirydowicz-nlawa}.
The triangular matrix $T_{k}^{(1)}$ is an approximation of the matrix
$(Q_{k}^TQ_{k})^{-1}$  
and two iterations result in a $T_{k}^{(2)}$ that is almost symmetric.
Note that here, the matrix $Q_k$ refers to the $A = QR$ factorization via
Gram-Schmidt, whereas $V_k$ refers to the orthogonal matrix produced in
the Arnoldi-$QR$ expansion.
Giraud et al.~\cite{gigl:simax:04} demonstrated how a 
rank-$k$ correction could be applied in a post-processing step 
to recover orthogonality by computing the
polar decomposition of ${Q}_{k-1}$, 
the matrix exhibited by Bj\"orck and Paige \cite{Bjorck92}.

Our paper is organized as follows. 
Low-synchronization Gram-Schmidt algorithms for the $QR$ factorization
are reviewed in Section 2 and multiplicative 
iterations are applied to solve the normal equations.
A rounding error analysis of the iterated Gauss-Seidel 
Gram-Schmidt
algorithm is presented in Section 3, leading to bounds on the
error and the orthogonality of the columns of $V_{k-1}$. In Section 4, 
we present the iterated Gauss Seidel GMRES algorithm and prove its backward stability 
and also derive a variant that requires only a single synchronization. 
Further, the relationship with Henrici's departure from normality is explored. 
Finally, numerical experiments on challenging problems studied
over the past thirty-five years are presented in Section 5.

{\it Notation.} Lowercase bold letters denote column vectors 
and uppercase letters are matrices (e.g., $\vect{v}$ and $A$,
respectively).  We use $A_{ij}$ to represent the $(i,j)$ scalar entry 
of a matrix $A$, and $\vect{a}_k$ denotes the $k$--th column of $A$.
$A_k$ is a block partition up to the $k$--th column of a matrix.
Subscripts indicate the approximate solution and corresponding residual 
(e.g., $\vect{x}_k$ and $\vect{r}_k$) 
of an iterative method at step $k$. Throughout this
article, the notation $U_k$ (or $L_k$) and $U_s$ 
(or $L_s$) will explicitly refer to {\it strictly upper/lower triangular
matrices}.\footnote{We note that the distinction between these 
two notations is crucial.  For $U_k$, the size of the strictly upper
triangular matrix changes with $k$, whereas the size of $U_s$ remains
fixed.}  Vector notation indicates a subset of the rows
and/or columns of a matrix; e.g. $V_{1:m+1,1:k}$ denotes the 
first $m + 1$ rows and $k$ columns of the matrix $V$ and the
notation $V_{:,1:k}$ represents the first $k$ columns 
of $V$.  $H_{k+1,k}$ represents an $(k+1)\times k$ upper Hessenberg matrix, and 
$h_{k+1,k}$ refers to a matrix element.  
Bars denote computed quantities (e.g., $\bar{Q}_{k-1}$ or $\bar{R}_{k-1}$). 
In cases where standard notation in 
the literature is respected that may otherwise conflict with the 
aforementioned notation, this will be explicitly indicated.
Note that the residual vector at iteration $k$ is denoted $\vect{r}_k$,
whereas $\vect{r}_{:,k}$ is the $k$--th column of the matrix $R$.

\section{Low-synchronization Gram-Schmidt Algorithms}

Krylov subspace methods for solving linear systems are often required for extreme-scale
applications on parallel machines with many-core accelerators. 
Their strong-scaling is limited by the number and frequency of
global reductions in the form of {\tt MPI\_Allreduce} operations and these communication
patterns are expensive \cite{Lockhart2022}.  
Low-synchronization algorithms are based on the ideas in 
Ruhe \cite{Ruhe83}, and are designed such that they require only 
one reduction per iteration to normalize each vector and apply projections. 
The Gram-Schmidt projector applied to $\vect{a}_k$, the $k$-th column of the $n \times m$ matrix $A$
in the factorization $A = QR$, can be written as 
\[
P_{k-1}\vect{a}_k = \vect{a}_k- Q_{k-1}\vect{r}_{1:k-1,k} = \vect{a}_k - 
Q_{k-1}\:(Q_{k-1}^T Q_{k-1} )^{-1}\:Q_{k-1}^T\vect{a}_k,
\]
where the vector $\vect{r}_{1:k-1,k}$ is the solution of 
the normal equations 
\begin{equation}
Q_{k-1}^T Q_{k-1}\vect{r}_{1:k-1,k} = Q_{k-1}^T\vect{a}_k.
 \label{eq:normal}
\end{equation}
Ruhe \cite{Ruhe83} established that the
MGS algorithm employs a {\it multiplicative}
Gauss-Seidel relaxation scheme with matrix splitting 
$Q_{k-1}^TQ_{k-1} = M_{k-1} - N_{k-1}$, 
where $M_{k-1}= I + L_{k-1}$ and $N_{k-1} = -L_{k-1}^T$. 
The iterated CGS-2
is an {\it additive} Jacobi relaxation.

The inverse compact $WY$ form  was derived in 
\'{S}wirydowicz et al.~\cite{2020-swirydowicz-nlawa}, with
strictly lower triangular matrix $L_{k-1}$. 
Specifically, these inverse compact $WY$ algorithms batch 
the inner-products together and compute the last row of $L_{k-1}$ as
\begin{equation}\label{eq:katvec}
L_{k-1,1:k-2} = \vect{q}_{k-1}^T\:Q_{k-2}.  
\end{equation}
The approximate projector $P_{k-1}^{(1)}$ and correction matrix are given by
\begin{equation}
P_{k-1}^{(1)} =  I - Q_{k-1}\:T_{k-1}^{(1)}\:Q_{k-1}^T, \quad 
T_{k-1}^{(1)} = (I + L_{k-1})^{-1}= M_{k-1}^{-1},
\label{eq:correct}
\end{equation}
respectively, and correspond to one iteration for 
the normal equations (\ref{eq:normal}) with the zero vector
as the initial guess. 
Iteration $k$ of the Gram-Schmidt algorithm with two Gauss-Seidel iterations is 
given as Algorithm \ref{fig:PSGSA} .
The compact and inverse compact $WY$ forms of the correction
matrix appearing in \eqref{eq:correct} 
can be derived as follows. The projector can be expressed as
\begin{eqnarray*}
P_{k-1}^{(1)} & = & (\: I - \vect{q}_{k-1}\vect{q}_{k-1}^T \:)\:(\: I - Q_{k-2}\:T_{k-2}^{(1)}\: Q_{k-2}^T\:) \\
  & = & I - \vect{q}_{k-1}\vect{q}_{k-1}^T - Q_{k-2}\:T_{k-2}^{(1)} \:Q_{k-2}^T 
  + \vect{q}_{k-1}\vect{q}_{k-1}^T\:Q_{k-2}\:T_{k-2}^{(1)} \:Q_{k-2}^T
\end{eqnarray*}
where $T_{k-1}^{(1)}$ is the lower triangular matrix of basis vector inner products.
The $WY$ representation is then given by the matrix form
\begin{equation}
P_{k-1}^{(1)} = I -
\left[ \begin{array}{cc} Q_{k-2} & \vect{q}_{k-1} \end{array} \right]
\left[
\begin{array}{cc}
T_{k-2}^{(1)}     & 0 \\
-\vect{q}_{k-1}^T\:Q_{k-2}\:T_{k-2}^{(1)} & 1
\end{array}
\right]
\left[ \begin{array}{c} Q_{k-2}^T \\ \vect{q}_{k-1}^T \end{array} \right]
\label{eq:WY}
\end{equation}
For the inverse compact $WY$ MGS, given the correction matrix $T_{k-1}^{(1)}$
generated by the recursion
\begin{equation}
(T^{(1)}_{k-1})^{-1} = I + L_{k-1} = 
\left[
\begin{array}{cc}
(T^{(1)}_{k-2})^{-1}      & 0 \\
\vect{q}_{k-1}^TQ_{k-2}  & 1
\end{array}
\right],
\end{equation}
it is possible to prove that
\begin{equation}
\left[
\begin{array}{cc}
(T^{(1)})_{k-2}^{-1}      & 0 \\
\vect{q}_{k-1}^TQ_{k-2}  & 1
\end{array}
\right]
\left[
\begin{array}{cc}
T^{(1)}_{k-2}                   & 0 \\
-\vect{q}_{k-1}^TQ_{k-2}T^{(1)}_{k-2}  & 1
\end{array}
\right] =
\left[
\begin{array}{cc}
I_{k-2} & 0 \\
0 & 1
\end{array}
\right].
\label{eq:block}
\end{equation}
Therefore, the lower triangular correction matrix 
$T^{(1)}_{k-1} = (I+L_{k-1})^{-1}$ from
\'{S}wirydowicz et al.~\cite{2020-swirydowicz-nlawa}
is equivalent to the compact $WY$ matrix $L_1$ from
Bj\"orck \cite[page 313]{Bjorck94}.

To prove that the correction matrix $T_{k-1}^{(2)}$ is close to 
a symmetric matrix, consider the matrix
\[
T_{k-1}^{(2)}  =  ( I + L_{k-1})^{-1} - 
(I + L_{k-1})^{-1} L_{k-1}^T (I + L_{k-1})^{-1} =
M_{k-1}^{-1} \left[ I + N_{k-1} M_{k-1}^{-1} \right].
\]
From this, we have $T_{k-1}^{(2)} =  T_{k-1}^{(1)} - T_{k-1}^{(1)} L_{k-1}^T T_{k-1}^{(1)}$.
From the block inverse \eqref{eq:block}, 
\[
T_{k-1}^{(1)} = 
\left[
\begin{array}{cc}
  T_{k-2}^{(1)}                      &  0 \\
 -\vect{q}_{k-1}^T Q_{k-2} T_{k-2}^{(1)} & 1
\end{array}
\right], \quad
L_{k-1}^T =
\left[
\begin{array}{cc}
      L_{k-2}^T    &  Q_{k-2}^T \vect{q}_{k-1}  \\
      0    & 0
\end{array}
\right].
\]
By subtracting the matrices above, and dropping ${\cal O}(\eps^2)$ terms,
the block symmetric matrix 
\begin{equation}
T_{k-1}^{(2)}  \approx
\left[
\begin{array}{cc}
  T_{k-2}^{(2)}              &  -Q_{k-2}^T \:\vect{q}_{k-1}  \\
 -\vect{q}_{k-1}^T\: Q_{k-2}   & 1
\end{array}
\right]
\label{eq:T2}
\end{equation}
is obtained and its singular values and eigenvalues remain close to one. However, as $k$ increases
$\sigma_{\max}(T_{k-1}^{(1)}) > 1$ and $\sigma_{\min}(T_{k-1}^{(1)}) < 1$, and
additionally $T_{k-1}^{(1)}$ is non-normal.
We provide definitions of the departure from normality in Section 4. 

 Algorithm \ref{fig:PSGSA}  corresponds 
to MGS-2 (two passes of MGS) and  is equivalent to two Gauss-Seidel iterations
in exact arithmetic; see Ruhe \cite{Ruhe83}. 
Following the MATLAB notation, the algorithm at the $k$-th step generates $QR$ factorization 
$A_{:, 1:k} = Q_{:,1:k}\:R_{1:k, 1:k}$, for $k =1, \dots$,
where $A_{:, 1:k} = [\:A_{:, k-1}, \: \vect{a}_k\:]$. 
The  matrices thus have dimensions $A_{:,1:k}$: $n\times k$, 
$Q_{:, 1:k}$: $n\times k$, and ${R}_{1:k,1:k}$: $k \times k$.
The diagonal elements are $\gamma_{k-1} = R_{k-1,k-1}$.
Two initial steps prime a depth-two pipeline.
\begin{algorithm}[H]
\centering
\begin{algorithmic}[1]
\State{$\vect{w}_1 := \vect{a}_1$, $R_{1,1} = \|\vect{w}_1\|_2$, 
$\vect{q}_1 := \vect{w}_1 / R_{1,1}$, $Q_1 = \vect{q}_1 $ \label{21:01}} 
\State{$R_{1,2} =  \:\vect{q}_{1}^T \:  \vect{a}_2$, 
$\vect{w}_2 := \vect{a}_2 - R_{1,2} \: \vect{q}_1 $, $Q_{2} = [\:Q_{1}, \: \vect{w}_2\:]$ 
\label{21:02} }
\For{$k = 3, \dots, m$ \label{21:03}  } 
     \Statex
    \State{$[\: L_{1:k-2,k-1}^T ,\: \vect{r}_{1:k-1,k}^{(0)}] =  
    [\:Q_{k-2}^T \: \vect{w}_{k-1}\ \:,\: Q_{k-1}^T \: \vect{a}_k\: ]$
        \Comment{Global Synchronization} \label{21:04} } 
  \State{$\gamma_{k-1} = \|\vect{w}_{k-1}\|_2$ \label{21:05} }
       \Statex
     \State{$\vect{q}_{k-1} = \vect{w}_{k-1} \: / \: \gamma_{k-1}$ 
     \Comment{Lagged Normalization} \label{21:06} } 
     \State{$\vect{r}_{1:k-1,k}^{(0)} = \vect{r}_{1:k-1,k}^{(0)} \: /  \: \gamma_{k-1}$ \label{21:07} } 
     \State{$L_{k-1,1:k-2} = L_{k-1,1:k-2} \: / \: \gamma_{k-1}$ \label{21:08} }
     \Statex
    \State{$\vect{r}_{1:k-1,k}^{(1)} = 
     (\:I + L_{k-1}\:)^{-1}\:\vect{r}_{1:k-1,k}^{(0)}$ \label{21:09}  }
     \State{$\vect{u}_k  = \vect{a}_k - {Q}_{k-1}\:\vect{r}_{1:k-1,k}^{(1)}$ 
     \Comment{ First Gauss-Seidel} \label{21:10} }
    \Statex
    \State{$\vect{r}_{1:k-1,k}^{(2)} =  \:Q_{k-1}^T \:  \vect{u}_k $
        \Comment{Global Synchronization} \label{21:11} }
    \State{$\vect{r}_{1:k-1,k}^{(3)} = 
     (\:I + L_{k-1}\:)^{-1}\: \vect{r}_{1:k-1,k}^{(2)}$ \label{21:12} }
        \State{$\vect{w}_{k} = \vect{u}_{k} 
        - Q_{k-1}\:\vect{r}_{1:k-1,k}^{(3)}$
     \Comment{ Second Gauss-Seidel } \label{21:13} }
     \Statex
     \State{$R_{1:k-1,k} = \vect{r}^{(1)}_{1:k-1,k} + \vect{r}^{(3)}_{1:k-1,k}$ \label{21:14}} 
     \State{$Q_{k} = [\:Q_{k-1}, \: \vect{w}_k\:]$ \label{21:15} }
     \Statex
     \EndFor
\end{algorithmic}
\caption{\label{fig:PSGSA} Inverse Compact $WY$ Two-Reduce 
Gauss-Seidel Gram-Schmidt}
\end{algorithm}
The projection is computed after one 
iteration, where $\vect{r}_{1:k-1,k}^{(1)}$ is an approximate
solution to the  normal equations, and can be expressed as 
\begin{equation}
\vect{u}_k  = \vect{a}_k - {Q}_{k-1}\:\vect{r}_{1:k-1,k}^{(1)} = 
\vect{a}_k -{Q}_{k-1} \: M_{k-1}^{-1}\:\vect{r}_{1:k-1,k}^{(0)}.
\label{eq:proj}
\end{equation}
After two iterations, the projection is given by
\begin{equation}
\vect{w}_k = 
\vect{a}_k - Q_{k-1} \vect{r}_{1:k-1,k}^{(1)}
- Q_{k-1} M_{k-1}^{-1} N_{k-1} \vect{r}_{1:k-1,k}^{(1)},
\label{eq:proj1}
\end{equation}
where 
$R_{k,k} = \| \vect{w}_k \|_2$ and $\vect{q}_{k} = \vect{w}_k  / R_{k,k}$.
The corresponding projector and correction matrix for two iterations are then found to be
\[
P_{k-1}^{(2)} =  I - Q_{k-1} T_{k-1}^{(2)} Q_{k-1}^T, \quad
T_{k-1}^{(2)} = (I+ L_{k-1})^{-1}  [I - L_{k-1}^T (I+L_{k-1})^{-1}]  
= M_{k-1}^{-1} [  I + N_{k-1} M_{k-1}^{-1} ],
\]
and the correction matrix is close to a symmetric matrix. 
The projection is applied across Steps \ref{21:09} through \ref{21:13}
in Algorithm \ref{fig:PSGSA}, and
employs the following recurrence:
\begin{eqnarray}
\vect{{r}}_{1:k-1,k}^{(0)} & = & Q_{k-1}^T \vect{a}_k, 
\nonumber \\
\vect{{r}}_{1:k-1,k}^{(1)}  & = & 
M_{k-1}^{-1}\:\vect{{r}}_{1:k-1,k}^{(0)},
\label{eq:recur}\\
\vect{{r}}_{1:k-1,k}^{(1)} + \vect{{r}}_{1:k-1,k}^{(3)} & = & 
\left[ I +
M_{k-1}^{-1} N_{k-1} \right]
\vect{{r}}_{1:k-1,k}^{(1)}.
\nonumber
\end{eqnarray}
It follows from (\ref{eq:proj}) and (\ref{eq:proj1}) that the vectors
$ \vect{u}_k$ and $\vect{w}_k$ are related as 
\begin{equation}
 \vect{w}_k  =  \vect{u}_k  -
 Q_{k-1}  M_{k-1}^{-1}\:N_{k-1}\:\vect{{r}}_{1:k-1,k}^{(1)}.    
\end{equation}
The orthogonality of the vectors  $ \vect{u}_k$ and $\vect{w}_k$  
with respect to columns of $Q_{k-1}$ is given as
\begin{eqnarray}
Q_{k-1}^T \vect{u}_k & = & Q_{k-1}^T\: \vect{a}_k -
Q_{k-1}^T Q_{k-1} M_{k-1}^{-1} Q_{k-1}^T \vect{a}_{k}  
= N_{k-1} M_{k-1}^{-1} Q_{k-1}^T \vect{a}_{k},\\ 
Q_{k-1}^T \vect{w}_k & = & Q_{k-1}^T \vect{u}_k  -
Q_{k-1}^T Q_{k-1} M_{k-1}^{-1} N_{k-1} M_{k-1}^{-1}
Q_{k-1}^T \vect{a}_{k}=
(N_{k-1} M_{k-1}^{-1})^2 Q_{k-1}^T \vect{a}_{k}.
\end{eqnarray}

\section{Rounding error analysis}

We now present a rounding error analysis of the low-synchronization  
Gram-Schmidt algorithm with two Gauss-Seidel iterations. 
Bj\"orck \cite{Bjork_1967} and Bj\"orck and Paige \cite{Bjorck94}
analyzed the traditional MGS algorithm. Our approach is
more akin to the analysis of the CGS-2 algorithm 
presented by Giraud et al.~\cite{Giraud2005}.
We begin with the derivation of the formulas for the computed quantities in the $QR$ 
factorization by invoking an induction argument. 
In particular, the triangular solves with $M_{k-1} = I + L_{k-1}$ 
implied by Step \ref{21:09} and Step \ref{21:12} of Algorithm \ref{fig:PSGSA} 
are backward stable 
as shown by Higham \cite{Higham1989}, where 
\[
(M_{k-1} + E_{k-1})\:\vect{x} = \vect{b},\quad \|E_{k-1}\|_2 \le {\cal O}(\eps)\:\|M_{k-1}\|_2.
\]
The associated error terms and bounds for the recurrence relations 
in Algorithm \ref{fig:PSGSA}
are then given as
\begin{align*}
\vect{\bar{r}}_{1:k-1,k}^{(0)}  & =  \tq_{k-1}^T \vect{a}_k + \vect{e}_k^{(0)}, &
\|\vect{e}_k^{(0)}\|_2 & \le {\cal O}(\eps) \|\tq_{k-1}\|_2\:\|\vect{a}_k\|_2, \\
( M_{k-1} + E_{k-1}^{(1)} \:)\:\vect{\bar{r}}_{1:k-1,k}^{(1)}  &  = 
\vect{\bar{r}}_{1:k-1,k}^{(0)}, &
\|E_{k-1}^{(1)}\|_2 & \le {\cal O}(\eps) \|M_{k-1}\|_2, \\
\vect{\bar{r}}_{1:k-1,k}^{(2)}  & =    \tq_{k-1}^T \vect{\bar{u}}_k + \vect{e}_k^{(2)}, &
\|\vect{e}_k^{(2)}\|_2 & \le {\cal O}(\eps) \|\tq_{k-1}\|_2\:\|\vect{\bar{u}}_k\|_2, \\
( M_{k-1} + E_{k-1}^{(2)} \:)\:\vect{\bar{r}}_{1:k-1,k}^{(3)}   & = 
\vect{\bar{r}}_{1:k-1,k}^{(2)},  &
\|E_{k-1}^{(2)}\|_2 & \le {\cal O}(\eps) \|M_{k-1}\|_2. 
\end{align*}
Throughout the entire paper, we employ the ${\cal O}(\eps)$ notation, 
as in Paige and Strako\v{s} \cite{Paige2002}, which denotes
a small multiple of machine precision where the proportionality constant is a low
degree polynomial in the matrix dimension $n$, and the iteration $k \ll n$.
Our analysis is based on standard assumptions on the computer arithmetic and
these parameters as used in Higham \cite{HighamBook2002}.
The computed coefficient vector $\vect{\bar{r}}_{1:k-1,k}^{(1)}$ satisfies 
\begin{equation*}
\vect{\bar{r}}_{1:k-1,k}^{(1)}  =  ( M_{k-1} + E_{k-1}^{(1)})^{-1}\: \vect{\bar{r}}_{1:k-1,k}^{(0)}
\end{equation*}
and its 2-norm satisfies the  inequality
\begin{eqnarray*}
\|\vect{\bar{r}}_{1:k-1,k}^{(1)}\|_2 & \le & 
\:\| ( M_{k-1} + E_{k-1}^{(1)} )^{-1}\:
\|\:\tq_{k-1}^T \vect{a}_k + \vect{e}_k^{(0)} \:\|_2
\\
& \le &
\frac{ 1 + {\cal O}(\eps) }
{\sigma_{k-1}(M_{k-1}) - \|E_{k-1}^{(1)}\|_2 } \|\tq_{k-1}\|_2\:\|\vect{a}_k\|_2 \\
& \le &
\:\frac{(\: 1 + {\cal O}(\eps)\:)\:\|M_{k-1}^{-1}\|_2 }
{1- {\cal O}(\eps)\kappa(M_{k-1}) }\|\tq_{k-1}\|_2\:\|\vect{a}_k\|_2
\end{eqnarray*}
assuming that ${\cal O}(\eps)\kappa(M_{k-1})<1$. 
The computed intermediate coefficient vector $\vect{\bar{r}}_{1:k-1,k}^{(3)}$ satisfies  
the formula
\begin{equation*}
\label{eq:plus}
\vect{\bar{r}}_{1:k-1,k}^{(3)}   = 
(M_{k-1} + E_{k-1}^{(2)})^{-1} \:\vect{\bar{r}}_{1:k-1,k}^{(2)},  
\end{equation*}
%
and an upper bound is given by
\begin{eqnarray*}
\|\vect{\bar{r}}_{1:k-1,k}^{(3)}\|_2 & \le &
\| ( M_{k-1} + E_{k-1}^{(2)} )^{-1} \|_2   \:
\|\:\tq_{k-1}^T \vect{\bar{u}}_k + \vect{e}_k^{(2)}\:\|_2  \|_2 \\
& \le & 
\frac{(1 +  {\cal O}(\eps) )\:  \|M_{k-1}^{-1} \|_2}
{1- {\cal O}(\eps)\kappa(M_{k-1}) } \|\tq_{k-1}\|_2\:\|\vect{\bar{u}}_k\|_2.
\end{eqnarray*}

The computed vector $\vect{\bar{u}}_k$ in Step \ref{21:10} of Algorithm \ref{fig:PSGSA} satisfies 
\begin{eqnarray}
\vect{\bar{u}}_k & = & \vect{a}_k - \tq_{k-1} \: \vect{\bar{r}}_{1:k-1,k}^{(1)} +
\vect{e}_k^{(1)}, 
\label{eq:baru} \\
\|\vect{e}_k^{(1)}\|_2 & \le & {\cal O}(\eps) 
\left[ \:
\|\vect{a}_k\|_2 + \|\tq_{k-1}\|_2\:\|\vect{\bar{r}}_{1:k-1,k}^{(1)}\|_2 \: \right].
\label{eq:ek1}
\end{eqnarray}
Similarly, the computed form of the projection in Step \ref{21:13} of Algorithm \ref{fig:PSGSA} 
satisfies
\begin{eqnarray}
\vect{\bar{w}}_k & = &
 \vect{\bar{u}}_k - \tq_{k-1} \: \vect{\bar{r}}_{1:k-1,k}^{(3)}  
+   \vect{e}^{(3)}_k,
\label{eq:wbark} \\
\|\vect{e}_k^{(3)}\|_2 & \le & {\cal O}(\eps) 
\left[ \:
\|\vect{\bar{u}}_k\|_2 + \|\tq_{k-1}\|_2\:\|\vect{\bar{r}}_{1:k-1,k}^{(3)}\|_2 \: \right].
\label{eq:ek2}
\end{eqnarray}
Summarizing (\ref{eq:baru}) and (\ref{eq:wbark}),
the lower bound for the norm of $\vect{\bar{u}}_k$ is now 
determined from the recurrences for the computed quantities
\begin{eqnarray}
A_{k-1} & = & \tq_{k-1}\:
\left (\bar{R}^{(1)}_{k-1} + \bar{R}^{(3)}_{k-1} \right ) - 
F_{k-1}, \  F_{k-1} =  \left[\begin{array}{ccc}
\vect{e}^{(1)}_{1}+  \vect{e}^{(3)}_{1} &
\dots,  &  \vect{e}^{(1)}_{k-1}+  \vect{e}^{(3)}_{k-1} 
\end{array}\right],\label{eq:rep1} \\
\vect{a}_k & = & \tq_{k-1}\:\vect{\bar{r}}_{1:k-1,k}^{(1)} + 
\vect{\bar{u}}_k - \vect{e}^{(1)}_{k},
\label{eq:rep2} 
\end{eqnarray}
where $A_{k-1}$ is given as $A_{k-1} = [\vect{a}_1, \dots, \vect{a}_{k-1} ] $. 
Combining these into an augmented matrix form, we obtain
\begin{equation}
A_k = \left[\begin{array}{cc} A_{k-1}, & \vect{a}_k \end{array} \right] =
\tq_{k-1}
\left[\begin{array}{cc} \bar{R}^{(1)}_{k-1} + \bar{R}^{(3)}_{k-1} , & \vect{\bar{r}}_{1:k-1,k}^{(1)} 
\end{array} \right] - 
\left[
\begin{array}{cc} F_{k-1}, & \vect{e}^{(1)}_{k} - \vect{\bar{u}}_k  \end{array}
\right],
\label{eq:normA}
\end{equation}
which results in the bound 
\begin{equation}
\left\| \left[ \begin{array}{cc}  F_{k-1}, 
&   \vect{e}^{(1)}_{k}    \end{array} \right ]
 \right\|_2
+ \|\vect{\bar{u}}_k \|_2 \ge 
\left\| \left[\begin{array}{cc} F_{k-1}, &   \vect{e}^{(1)}_{k}  - \vect{\bar{u}}_k 
\end{array} \right] \right\|_2 \ge \sigma_{\min}(A_k).
\label{eq:upper}
\end{equation}
Thus, the 2-norm of the projected vector $\vect{\bar{u}}_k$ appearing
in \eqref{eq:wbark} is bounded below as
\begin{equation}
\| \vect{\bar{u}}_k \|_2 \ge \sigma_{\min}(A_k) - \left  \|\left[\begin{array}{cc} F_{k-1}, 
&   \vect{e}^{(1)}_{k} 
\end{array} \right] \right \|_2.
\label{eq:lower}
\end{equation}
Invoking the inequalities  (\ref{eq:baru}), (\ref{eq:wbark}) and the upper bounds on
$\|\vect{\bar{r}}_{1:k-1,k}^{(1)}\|_2$ and $\|\vect{\bar{r}}_{1:k-1,k}^{(3)}\|_2$, 
together with the bounds 
$\| \vect{\bar{u}}_k \|_2 \approx \| \vect{a}_k \|_2 \leq \| A_k\|_2$, 
it follows that 
\begin{equation}
\label{boundonf}
 \left  \|\left[\begin{array}{cc} F_{k-1}, &   \vect{e}^{(1)}_{k} 
\end{array} \right] \right \|_2
\leq {\cal O}(\eps) \frac{(\: 1 + {\cal O}(\eps)\:)\:\|M_{k-1}^{-1}\|_2 }
{(1- {\cal O}(\eps)\kappa(M_{k-1}) )^2}\|\tq_{k-1}\|^2_2\:\| A_k\|_2.
\end{equation}

In order to derive an upper bound for
$\|\tq_{k-1}^T\:\vect{\bar{w}}_{k}\:\|_2 \: / \: \|\vect{\bar{w}}_{k}\:\|_2$,
assume inductively that $\| N_{k-1} \|_2  \leq {\cal O}(\eps)$. 
Note that the term  ${\cal O}(\eps)$ depends on the iteration number $k$.
Then it is evident that $\|\tq_{k-1} \|^2 \leq 1+ {\cal O}(\eps) $ and 
$\| M_{k-1} \|_2 \leq 1 + {\cal O}(\eps)$, with $ \sigma_{\min}(M_{k-1}) 
\geq 1 - {\cal O}(\eps)$. Consequently, it follows that
\begin{equation*}
    \kappa(M_{k-1}) \leq \frac{1 + {\cal O}(\eps)}{1 - {\cal O}(\eps)}.
\end{equation*} 
The term $\tq_{k-1}^T\:\vect{\bar{u}}_{k}$ can be written in the form
\begin{eqnarray*}
 \tq_{k-1}^T \vect{\bar{u}}_{k}
 & = &  \tq_{k-1}^T\vect{a}_k - \tq_{k-1}^T\tq_{k-1} \vect{\bar{r}}_{1:k-1,k}^{(1)} -
 \tq_{k-1}^T \: \vect{e}_k^{(1)} \\ 
 & = & \tq_{k-1}^T \vect{a}_{k} - 
(M_{k-1} - N_{k-1}) ( M_{k-1} + E_{k-1}^{(1)})^{-1}\: 
\vect{\bar{r}}_{1:k-1,k}^{(0)} - \tq_{k-1}^T\:\vect{e}_k^{(1)}
\\
& = & 
(N_{k-1} + E_{k-1}^{(1)} ) ( M_{k-1} + E_{k-1}^{(1)} )^{-1}\:
\vect{\bar{r}}_{1:k-1,k}^{(0)} - \tq_{k-1}^T\:\vect{e}_k^{(1)} - \vect{e}_k^{(0)}.
\end{eqnarray*} 
Thus, the quotient  
$\|\:\tq_{k-1}^T\:\vect{\bar{u}}_{k}\:\|_2 \: / \: \|\:\vect{\bar{u}}_{k}\:\|_2$
can be bounded as 
\begin{equation}
\frac{\|\tq_{k-1}^T \vect{\bar{u}}_{k}\:\|_2}
{ \|\vect{\bar{u}}_{k}\:\|_2}  \leq
\|  N_{k-1} + E_{k-1}^{(1)} \|_2 \| (  M_{k-1} + E_{k-1}^{(1)} )^{-1} \|_2 
\frac{\| \vect{\bar{r}}_{1:k-1,k}^{(0)} \|_2}{\|\vect{\bar{u}}_{k}\:\|_2} + 
\frac{\| \tq_{k-1}^T \vect{e}_k^{(1)} \|_2}{\|\vect{\bar{u}}_{k}\:\|_2} 
+ \frac{\|\vect{e}_k^{(0)}\|_2}{\|\vect{\bar{u}}_{k}\:\|_2},
\label{eq:Qu}
\end{equation}
and therefore
\begin{equation}
\frac{\|\tq_{k-1}^T \vect{\bar{u}}_{k}\:\|_2}
{ \|\vect{\bar{u}}_{k}\:\|_2}  \leq \left[\:
\frac{(1 + {\cal O}(\eps))  \| M_{k-1}^{-1}\|_2 }{ 1 - {\cal O}(\eps)\kappa(M_{k-1})} \:\right]
\: \| N_{k-1} \|_2
\frac{\| \vect{\bar{r}}_{1:k-1,k}^{(0)} \|_2}{\|\vect{\bar{u}}_{k}\:\|_2} + 
\frac{\| \tq_{k-1}^T \vect{e}_k^{(1)} \|_2}{\|\vect{\bar{u}}_{k}\:\|_2} 
+ \frac{\|\vect{e}_k^{(0)}\|_2}{\|\vect{\bar{u}}_{k}\:\|_2}.
\label{eq:Qu2}
\end{equation}
The vector 2-norm $\| \vect{\bar{r}}_{1:k-1,k}^{(0)} \|_2$ 
can be bounded according to 
\begin{equation}
\| \vect{\bar{r}}_{1:k-1,k}^{(0)} \|_2 \leq \left( 1 + {\cal O}(\eps) \right) \:
\|\tq_{k-1} \|_2 \: \|\vect{a}_k\|_2,
\end{equation}
and together with (\ref{eq:lower}) and (\ref{boundonf}) 
\begin{equation*}
 \frac{\| \vect{\bar{r}}_{1:k-1,k}^{(0)} \|_2}{\|\vect{\bar{u}}_{k}\:\|_2}  
 \lesssim   \frac{\kappa(A_k)}{1- {\cal O}(\eps)\kappa(A_{k})}
\end{equation*}
under the assumption that ${\cal O}(\eps) \kappa(A_{k}) < 1$. 
After one Gauss-Seidel iteration, the bound for the 
loss of orthogonality is proportional to 
\begin{equation}
\label{onesteploo}
\frac{\|\tq_{k-1}^T\:\vect{\bar{u}}_{k}\:\|_2}
{ \|\vect{\bar{u}}_{k}\:\|_2}
\lesssim \frac{{\cal O}(\eps) \kappa(A_k) }{1- {\cal O}(\eps)\kappa(A_{k})}.
\end{equation}
This corresponds to the upper bound on the loss of orthogonality 
originally derived by Bj\"orck \cite{Bjork_1967} and further
refined by Bj\"orck and Paige \cite{Bjorck92} for the MGS algorithm. 

Given \eqref{eq:wbark}, the projected vector is expanded as
\begin{eqnarray}
\ \tq_{k-1}^T \vect{\bar{w}}_{k}
 & = &  \tq_{k-1}^T\vect{\bar{u}}_k - \tq_{k-1}^T\tq_{k-1}\: \vect{\bar{r}}_{1:k-1,k}^{(3)} -
 \tq_{k-1}^T \: \vect{e}_k^{(3)} 
 \nonumber \\ 
& = &  \tq_{k-1}^T\:\vect{\bar{u}}_{k} - 
(M_{k-1} - N_{k-1}) ( M_{k-1} + E_{k-1}^{(2)})^{-1}\: 
\vect{\bar{r}}_{1:k-1,k}^{(2)} - \tq_{k-1}^T\:\vect{e}_k^{(3)}
\nonumber
\\
& = & 
(N_{k-1} + E_{k-1}^{(2)} ) ( M_{k-1} + E_{k-1}^{(2)} )^{-1}\:
\vect{\bar{r}}_{1:k-1,k}^{(2)} - 
\tq_{k-1}^T\:\vect{e}_k^{(3)} - \vect{e}_k^{(2)}. 
\nonumber 
\end{eqnarray}
Therefore, the quotient 
$\|\tq_{k-1}^T\:\vect{\bar{w}}_{k}\:\|_2 \: / \: \|\vect{\bar{w}}_{k}\:\|_2$
is bounded as 
\begin{eqnarray}
\frac{\|\tq_{k-1}^T\vect{\bar{w}}_{k}\:\|_2}
{ \|\vect{\bar{w}}_{k}\:\|_2}  & \leq &
\| N_{k-1} + E_{k-1}^{(2)}   \|_2  \| ( M_{k-1} + E_{k-1}^{(2)} )^{-1} \|_2
\frac{\|\vect{\bar{r}}_{1:k-1,k}^{(2)}\|_2}{ \|\vect{\bar{w}}_{k}\:\|_2}
\nonumber \\
& + & 
 \frac{\|\tq_{k-1}^T\:\vect{e}_k^{(3)} \|_2}{\|\vect{\bar{w}}_{k}\:\|_2}
 +
\frac{\|\vect{e}_k^{(2)} \|_2}{\|\vect{\bar{w}}_{k}\:\|_2}\nonumber \\
 & \leq &
\left [ \| N_{k-1} + E_{k-1}^{(2)}   \|_2  
\| ( M_{k-1} + E_{k-1}^{(2)} )^{-1} \|_2
\frac{\|\vect{\bar{r}}_{1:k-1,k}^{(2)}\|_2}{ \|\vect{\bar{u}}_{k}\:\|_2}
\right . 
\nonumber \\
& + & \left .
 \frac{\|\tq_{k-1}^T\:\vect{e}_k^{(3)} \|_2}{\|\vect{\bar{u}}_{k}\:\|_2}
 +
\frac{\|\vect{e}_k^{(2)} \|_2}{\|\vect{\bar{u}}_{k}\:\|_2}
\right ] \times \frac{\| \vect{\bar{u}}_{k} \|_2}{\| \vect{\bar{w}}_{k} \|_2}.
\nonumber
\end{eqnarray}
The term $\| \vect{\bar{r}}_{1:k-1,k}^{(2)} \|_2$ 
satisfies the bound 
\begin{equation}
\| \vect{e}_k^{(2)} \|=
\| \vect{\bar{r}}_{1:k-1,k}^{(2)} - 
\tq_{k-1}^T\:\vect{\bar{u}}_k  \|_2 \leq  {\cal O}(\eps)  \|
\tq_{k-1} \|_2 \|\vect{\bar{u}}_k\|_2.   
\end{equation}
After two Gauss-Seidel iterations, the essential 
orthogonality relation is obtained as
\begin{eqnarray}
\frac{\|\tq_{k-1}^T \vect{\bar{w}}_{k}\:\|_2}
{ \|\vect{\bar{w}}_{k}\:\|_2}  & \leq & \left[\:
\frac{(1 + {\cal O}(\eps))  \| M_{k-1}^{-1}\|_2 }{ 1 - {\cal O}(\eps)\kappa(M_{k-1})} 
\: \| N_{k-1} \|_2
\frac{\|\tq_{k-1}^T\:\vect{\bar{u}}_k \|_2 + \|\vect{e}_k^{(2)}\|_2}{\|\vect{\bar{u}}_{k}\:\|_2} \right . 
\nonumber \\ & + & \left .
\frac{\| \tq_{k-1}^T \vect{e}_k^{(3)} \|_2}{\|\vect{\bar{u}}_{k}\:\|_2} 
+ \frac{\|\vect{e}_k^{(2)}\|_2}{\|\vect{\bar{u}}_{k}\:\|_2}\:\right]
\times \frac{\| \vect{\bar{u}}_{k} \|_2}{\| \vect{\bar{w}}_{k} \|_2}.
\label{eq:Qw2}
\end{eqnarray}
By considering (\ref{eq:wbark}), an upper bound  
for the ratio $\| \vect{\bar{u}}_{k} \|_2 / \| \vect{\bar{w}}_{k} \|_2 $ 
is obtained as 
\begin{equation}
\frac{\| \vect{\bar{w}}_{k} \|_2}{\| \vect{\bar{u}}_{k} \|_2}    
\geq \frac{\| \vect{\bar{u}}_{k} \|_2}{\| \vect{\bar{u}}_{k} \|_2} -
\|\tq_{k-1} \|_2  \frac{\| \vect{\bar{r}}_{1:k-1,k}^{(3)} \|_2}{\| \vect{\bar{u}}_{k} \|_2}  
- \frac{\| \vect{e}_k^{(3)} \|_2}{\| \vect{\bar{u}}_{k} \|_2},
\end{equation}
and substituting the upper bound for $\| \vect{\bar{r}}_{1:k-1,k}^{(3)} \|_2$ this
can be rewritten as 
\begin{equation}
\frac{\| \vect{\bar{u}}_{k} \|_2}{\| \vect{\bar{w}}_{k} \|_2} 
\leq \frac{1}{ 1 -  {\cal O}(\eps)\kappa(A_{k}) }.
\end{equation}
An analogous bound was derived in 
Giraud et al.~\cite[eq. 3.33]{Giraud2005} for the classical Gram-Schmidt
algorithm (CGS-2) with reorthogonalization.
Thus, assuming ${\cal O}(\eps)\kappa(A_{k})< 1$ and given Bj\"orck's 
result (\ref{onesteploo}),  it follows that the loss of 
orthogonality is maintained at the level of machine precision:
\[
\|\tq_{k-1}^T\:\vect{\bar{q}}_{k}\|_2 \lesssim {\cal O}(\eps).
\]

\section{Iterated Gauss-Seidel GMRES}\label{sec:lowsynch}

The GMRES algorithm of Saad and Schultz \cite{Saad86} has a 
thirty-five year history and various alternative formulations
of the basic algorithm have been proposed over that time frame. 
A comprehensive review of these is presented 
by Zou \cite{Zhou2021}.
In particular, pipelined, $s$-step, and block algorithms
are better able to hide latency in parallel 
implementations; see, e.g., Yamazaki et al.~
\cite{2020-yamazaki-proceedings-of-siam-pp20}. 
The low-synchronization MGS-GMRES algorithm
described in \'{S}wirydowicz et al.~\cite{2020-swirydowicz-nlawa}
improves parallel
strong-scaling by employing one global reduction for each GMRES 
iteration; see Lockhart et al.~\cite{Lockhart2022}.
A review of compact $WY$ Gram Schmidt algorithms and their computational
costs is given in~\cite{DCGS2Arnoldi}.
Block generalizations of the DGCS-2 and CGS-2 algorithm are presented in Carson
et al.~\cite{Carson2022, CarsonRoz2021}.  In \cite{CarsonRoz2021} the authors generalize
the Pythagorean theorem to block form and derive {\tt BCGS-PIO} and {\tt
BCGS-PIP} algorithms with the more favorable communication patterns described
herein. An analysis of the backward stability of the these block Gram-Schmidt
algorithms is also presented. Low-synch 
iterated Gauss-Seidel GMRES algorithms are now presented  with
one and two global reductions.

The low-synchronization DCGS-2 algorithm introduced by 
\'{S}wirydowicz \cite{2020-swirydowicz-nlawa}
was employed to compute the Arnoldi-$QR$ expansion in the Krylov-Schur 
eigenvalue algorithm by Bielich et al.~\cite{DCGS2Arnoldi}. 
The algorithm exhibits desirable
numerical characteristics including the computation of invariant subspaces
of maximum size for the Krylov-Schur algorithm of Stewart \cite{Stewart86}. 
In the case of iterated Gauss-Seidel GMRES, the backward error 
analysis derived for the two-reduce Gauss-Seidel Gram-Schmidt Algorithm \ref{fig:PSGSA} can be 
applied to the IGS-GMRES algorithm. In the case of the DCGS-2 algorithm, the symmetric correction matrix 
$T_{k-1}$ was derived in Appendix 1 of \cite{2020-swirydowicz-nlawa} and is given by 
$
T_{k-1} = I - L_{k-1} - L_{k-1}^T.
$
This correction matrix was employed in $s$-step 
and pipe\-lined MGS-GMRES. When the matrix
$T_{k-1}$ is split into $I - L_{k-1}$ and $L_{k-1}^T$ and 
applied across two iterations of the DCGS-2 algorithm, the resulting loss of 
orthogonality is empirically observed to be ${\cal O}(\eps)$. 
Indeed, it was 
noted in Bielich et al.~\cite{DCGS2Arnoldi} that two iterations
of classical Gram-Schmidt (CGS)
are needed to achieve vectors orthogonal to the level ${\cal O}(\eps)$.

\subsection*{Two-Reduce IGS-GMRES}\label{sec:twoReduce}

The iterated Gauss-Seidel GMRES algorithm (IGS-GMRES)   
presented in Algorithm \ref{lowsynch}  requires
computing $\vect{v}_{k}$ in Step 4 and a norm in Step 6,
followed by vector scaling in Step \ref{41:07}. 
Two Gauss-Seidel iterations are applied in Step \ref{41:13} and \ref{41:16}.
The normalization for the Krylov vector $\vect{v}_{k}$ 
at iteration $k$ represents
the delayed scaling of the vector $\vect{w}_{k-1}$ in the
matrix-vector product $\vect{v}_{k} = A\vect{w}_{k-1}$. Therefore,
an additional Step \ref{41:08} is required, together with Step \ref{41:09}.
The subdiagonal element $\gamma_{k-1}$ 
in the Arnoldi-$QR$ expansion (\ref{eq:arnoldi})  is computed
in Step \ref{41:06} and the remaining entries $H_{1:k-1,k-1}$ are computed after the
second Gauss-Seidel iteration in Step 16 of Algorithm \ref{lowsynch}.
In order to lag the normalization to the next iteration, the
norm $\|\vect{w}_{k-1}\|_2$ is included in the global
reduction in Steps \ref{41:05} and \ref{41:06}.
Using MATLAB notation, 
the algorithm denoted as IGS-GMRES
at the $k$-th iteration step  computes the Arnoldi expansion 
$A\:V_{:, 1:k-1} = V_{:, 1:k}\:H_{1:k,1:k-1}$. 
This expansion holds only for $k < n-1$, as the $n$-th column of $V$ has not 
been normalized and the $n$-th column of $H$ has not been set.
The normalization coefficients are $\gamma_{k-1} = H_{k-1,k-2}$.
A lagged normalization leads to a pipeline of depth two. 
Thus, two initial iteration steps prime the pipeline.

\begin{algorithm}[H]
\centering
\caption{\label{lowsynch} Low-Synchronization Two-Reduce Iterated Gauss-Seidel GMRES }
\begin{algorithmic}[1]
\State{$\vect{r}_0 = b - A\:\vect{x}_0$, $\rho = \|\vect{r}_0 \|_2$, 
$\vect{w}_1 = \vect{r}_0$,
$\vect{v}_1 = \vect{w}_1 / \rho$ \label{41:01} }
\State{$\vect{w}_2 = A\vect{v}_1$, $H_{1,1} =  \:\vect{v}_{1}^T \: \vect{w}_2$, $\vect{w}_2 = \vect{w}_2 - H_{1,1} \: \vect{v}_1$, $V_{2} = [\:V_{1}, \: \vect{w}_2\:]$ \label{41:02} }
    \For{$k=3, \dots m$ \label{41:03} }
    \Statex
     \State{$\vect{v}_{k} = A\vect{w}_{k-1}$ \label{41:04} }
     \State{$[\: L_{1:k-2,k-1}^T,\: \vect{r}_{1:k-1,k}^{(0)} \:] =  
     \left[\: V_{k-2}^T\:\vect{w}_{k-1}, \:V_{k-1}^T\:\vect{v}_{k}\:\right]$ 
      \Comment{Global Synchronization} \label{41:05} } 
     \State{$\gamma_{k-1} = \|\vect{w}_{k-1}\|_2$ \label{41:06} }
     \State{$\vect{v}_{k-1} = \vect{w}_{k-1} \: / \: \gamma_{k-1}$ \label{41:07} }
     \Statex
     \State{$\vect{r}_{1:k-1,k}^{(0)} = \vect{r}_{1:k-1,k}^{(0)} \: / \: \gamma_{k-1}$
     \Comment{Scale for Arnoldi} \label{41:08} }
     \State{$\vect{r}_{k-1,k}^{(0)} = \vect{r}_{k-1,k}^{(0)} \: / \: \gamma_{k-1}$ \label{41:09} }
     \State{$L_{k-1,1:k-2} = L_{k-1,1:k-2} \: / \: \gamma_{k-1}$ \label{41:10} }
     \State{$\vect{v}_{k} = \vect{v}_{k} \: / \: \gamma_{k-1}$ \label{41:11} }
    \Statex
       \State{$\vect{r}_{1:k-1,k}^{(1)} = 
     (I + L_{k-1})^{-1}\:\vect{r}_{1:k-1,k}^{(0)}$ \label{41:12} }
     \State{$\vect{u}_k  = \vect{v}_{k} - {V}_{k-1}\:\vect{r}_{1:k-1,k}^{(1)}$ 
     \Comment{ First Gauss-Seidel} \label{41:13} }
     \Statex
     \State{$\vect{r}_{1:k-1,k}^{(2)} =  \:V_{k-1}^T \:  \vect{u}_k$
        \Comment{Global Synchronization} \label{41:14} }
    \State{$\vect{r}_{1:k-1,k}^{(3)} = 
      (I + L_{k-1})^{-1}\:\vect{r}_{1:k-1,k}^{(2)} $ \label{41:15} } 
     \State{$\vect{w}_{k} = \vect{u}_{k}
     - V_{k-1}\:\vect{r}_{1:k-1,k}^{(3)}$
     \Comment{ Second Gauss-Seidel} \label{41:16} }
     \Statex
     \State{$H_{1:k-1,k-1} = \vect{r}^{(1)}_{1:k-1,k} + \vect{r}^{(3)}_{1:k-1,k}$ \label{41:17} }
     \State{$V_{k} = [\:V_{k-1}, \: \vect{w}_k\:]$ \label{41:18} }
  \Statex
  \EndFor
   \State{$\vect{y}_m = {\rm argmin}  \:\|\: H_{m}\vect{y}_m - \rho\:\vect{e}_1 \:\|_2$ ,
 $\vect{x}_m = \vect{x}_0 + V_m\vect{y}_m$}
\end{algorithmic}
\end{algorithm}
For the least squares solution, we solve
$\vect{y}_m = {\rm argmin}  \:\|\: H_{m}\vect{y}_m - \rho\:\vect{e}_1 \:\|_2$ 
and then compute $\vect{x}_m = \vect{x}_0 + V_m\vect{y}_m$.
The subdiagonal element $\gamma_{k-1}$ is computed as 
$\|\vect{w}_{k-1}\|_2$ in Algorithm \ref{lowsynch}.  
To reduce the number of global synchronizations
the recurrence 
\begin{equation}
    \|\vect{w}_{k-1}\|_2^2 = 
    \|\: \vect{u}_{k-1}\:\|_2^2 - \|\:\vect{r}_{1:k-2,k-1}^{(2)}\:\|_2^2
\end{equation}
 would have to be employed. 
Therefore, we have
\[
\|\vect{w}_{k-1}\|_2 \ge \|\vect{u}_{k-1}\|_2 - \|\vect{r}^{(2)}_{1:k-2,k-1}\|_2.
\]
This means that the magnitude of the vector $\vect{w}_{k-1}$ is at least as large 
as the difference between the magnitudes of the vectors $\vect{u}_{k-1}$ and 
$\vect{r}^{(2)}_{1:k-2,k-1}$. Because $\| \vect{r}^{(2)}_{1:k-2,k-1}\|_2$ is small, 
it follows that $\|\vect{w}_{k-1}\|_2$ is approximately equal to $\|\vect{u}_{k-1}\|_2$.
Thus, if $\|\vect{r}^{(2)}_{1:k-2,k-1}\|_2$ is small, then the inequality is a tight bound, 
and $\|\vect{w}_{k-1}\|_2$ is approximately equal to $\|\vect{u}_{k-1}\|_2$.
Therefore, we can say that 
$\|Q_{k-1}^T\vect{q}_{k}\|_2 =\||Q_{k-1}^T\vect{w}_{k}\|_2 / \|\vect{w}_{k}\|_2$
is ${\cal O}(\eps)$.

%

\subsection*{Hybrid MGS-CGS GMRES}\label{sec:oneReduce}
We now derive a new variant of Algorithm \ref{lowsynch} that requires only a 
single synchronization point.
The vector $\vect{w}_{k-1}$ is available at Step 4 of Algorithm \ref{lowsynch}.
Alternatively, the vector $\vect{v}_k = A\vect{u}_{k-1}$ 
can be computed during the Gauss-Seidel iteration in Step 13,
by replacing Step 11 in Algorithm \ref{lowsynch}, where
$\vect{v}_{k} = A \vect{w}_{k-1} \: / \:\gamma_{k-1}$, and
$\vect{w}_{k-1} = \vect{u}_{k-1} - V_{k-2}\:V_{k-2}^T\:\vect{u}_{k-1}$,  i.e.,
\begin{equation}
\vect{v}_k  =  
\frac{1}{\gamma_{k-1}} \left[ \:
A \vect{u}_{k-1} - V_{k-1} H_{1:k-1,1:k-2}
\:\vect{r}^{(2)}_{1:k-2,k-1}
\:\right].
\label{eq:Aw} 
\end{equation}
The recurrence for Step 11 is instead written as 
$\vect{u}_k = A\:\vect{v}_{k-1} - V_{k-1} V_{k-1}^T  A\:\vect{v}_{k-1}$. It follows that
\begin{eqnarray*}
\vect{u}_{k}  & = &  \frac{1}{\gamma_{k-1}} \left[ \:
A \:\vect{u}_{k-1} - V_{k-2} \:  V_{k-2}^T\: A\: \vect{u}_{k-1} 
\:\right] - \frac{1}{\gamma_{k-1}} 
\vect{v}_{k-1}\:\vect{v}_{k-1}^T\:A\:\vect{u}_{k-1}
\nonumber \\
& + &  \frac{1}{\gamma_{k-1}} \: V_{k-1}\: \left( \:
I - V_{k-1}^T\:V_{k-1}\:\right) \:H_{1:k-1,1:k-2} \:\vect{r}^{(3)}_{1:k-2,k-1}.
\end{eqnarray*}
 Noting that 
$\vect{r}^{(0)}_{1:k-2,k} =  {V}_{k-2}^T\: A\: \vect{u}_{k-1}$,
it is possible to compute $\vect{r}^{(1)}_{k-1,k}$ as follows:
\begin{eqnarray*}
\vect{v}_{k-1}^T\:A\:\vect{u}_{k-1} & = &
\frac{1}{\gamma_{k-1}} \: \left[ \:
\vect{u}_{k-1} - V_{k-2}\:\vect{r}^{(2)}_{1:k-2,k-1} \:
\:\right]^T\:A\:\vect{u}_{k-1} \\
& = & 
\frac{1}{\gamma_{k-1}} \: \left[ \:
\vect{r}^{(0)}_{k-1,k} - L_{k-1,1:k-2}\:\vect{r}^{(0)}_{1:k-2,k}
\:\right],
\end{eqnarray*}
where $L_{k-1,1:k-2} = \vect{u}_{k-1}^T\:V_{k-2}$. This is ``Stephen's trick'' 
from \'{S}wirydowicz et al.~\cite[eq. 4]{2020-swirydowicz-nlawa} and
Bielich et al.~\cite{DCGS2Arnoldi} applied to the Arnoldi-$QR$ algorithm. 
The first projection step is then applied in Steps 9 and 11 using the compact $WY$
form given in equation (\ref{eq:WY}), i.e.,
\begin{equation}
\vect{u}_{k} =  A\:\vect{u}_{k-1} \:/\: \gamma_{k-1} - 
\left[ \begin{array}{cc} V_{k-2} & \vect{v}_{k-1} \end{array} \right]
\left[ \begin{array}{cc} I & 0  \\
-\gamma_{k-1}^{-1}\:L_{k-1,1:k-2} & 1 \end{array} \right]
\left[ \begin{array}{c} \vect{r}^{(0)}_{1:k-2,k} \\  
\vect{r}^{(0)}_{k-1,k}\end{array} \right],
\end{equation}
where the implied triangular inverse simplifies 
to $(\:I + L_{k-1}\:)^{-1} = I - L_{k-1}$, only when the last row
contains non-zero off-diagonal elements, as in Step 12 of Algorithm \ref{onereduce}. 
The correction matrix then takes the simplified form
\[
T_{k-1}^{(1)} =
\left[
\begin{array}{cc}
I    & 0 \\
-\vect{v}_{k-1}^T\:V_{k-2} & 1
\end{array}
\right],
\]
where the spectral radius $\rho_{k-1}$ of the matrix $M_{k-1}^{-1}\:N_{k-1}$ is 
identical to that of the two-reduce algorithm.
The projection matrix is given by
\begin{eqnarray*}
P_{k-1}^{(1)} & = & (\: I - \vect{v}_{k-1}\vect{v}_{k-1}^T \:)\:(\: I - V_{k-2}\:T_{k-1}^{(1)}\:V_{k-2}^T\:) \\
  & = & I - \vect{v}_{k-1}\vect{v}_{k-1}^T - V_{k-2} \:V_{k-2}^T 
  + \vect{v}_{k-1}\vect{v}_{k-1}^T\:V_{k-2} \:V_{k-2}^T .
\end{eqnarray*}
After substitution of this expression, it follows that
\[
\vect{u}_k  = \frac{1}{\gamma_{k-1}} \:\left[\: A\vect{u}_{k-1} - 
{V}_{k-1}\:\vect{r}^{(1)}_{1:k-1,k}
\:\right]. 
\]
The one-synch hybrid MGS-CGS GMRES
algorithm  is presented in Algorithm \ref{onereduce}.
The algorithm can be characterized by an MGS iteration at Step \ref{42:13},
combined with a lagged CGS iteration at Step \ref{42:10}.
The update of the Hessenberg matrix in Step \ref{41:14}
is the same as earlier. The Arnoldi expansion holds only for $k < n-1$.
We note that our backward stability analysis in the following section only applied to IGS-GMRES, Algorithm
\ref{lowsynch}, with two global synchronizations. A complete stability analysis of Algorithm \ref{onereduce} 
remains future work. 
\begin{algorithm}[htb]
\begin{algorithmic}[1]
\State{$\vect{r}_0 = b - A\:\vect{x}_0$, $\rho = \|\vect{r}_0 \|_2$, 
$\vect{w}_1 = \vect{r}_0$,
$\vect{v}_1 = \vect{w}_1 / \rho$ \label{42:01}} 
\State{$\vect{w}_2 = A\vect{v}_1$, $H_{1,1} =  \:\vect{v}_{1}^T \: \vect{w}_2$, 
$\vect{u}_2 = ( \vect{u}_2 - \vect{v}_1\: H_{1,1}\:)\:/\:H_{1,1}$, $V_{2} = [\:V_{1}, \: \vect{w}_2\:]$ \label{42:02} }
    \For{$k=3, \dots m$  \label{32}}
\Statex 
\State{$
\left [ \begin{array}{cc} \vect{r}^{(2)}_{1:k-2,k-1} & \vect{r}^{(0)}_{1:k-1,k}  \\ 
\|\vect{u}_{k-1}\:\|_2^2 &  \vect{r}^{(0)}_{k-1,k} \end{array} \right ] =
[\:{V}_{k-2},\: \vect{u}_{k-1}\:]^T\:
[\: \vect{u}_{k-1},\:A\vect{u}_{k-1}\:]$ \:\Comment{Global Synchronization} \label{42:04} }
\Statex
\State{$\vect{r}^{(3)}_{1:k-2,k-1} =
(\:I + L_{k-2}\:)^{-1}\:\vect{r}^{(2)}_{1:k-2,k-1}$ 
\Comment{Normal equations} \label{42:05} } 
\State{$\gamma_{k-1} = \left\{
(\: \|\vect{u}_{k-1}\:\|_2^2 - \|\:\vect{r}_{1:k-2,k-1}^{(2)}\:\|_2^2 \:) \right\}^{1/2}$ 
\Comment{lagged norm} \label{42:06} }
\Statex
\State{$\vect{r}_{1:k-1,k}^{(0)} = \vect{r}_{1:k-1,k}^{(0)} \: / \: \gamma_{k-1}$
\Comment{Scale for Arnoldi} \label{42:07}}
\State{$\vect{r}_{k-1,k}^{(0)} = \vect{r}_{k-1,k}^{(0)} \: / \: \gamma_{k-1}$ \label{42:08} }
\State{$L_{k-1,1:k-2} =  \vect{r}^{(2)}_{1:k-2,k-1} / \: \gamma_{k-1}$ 
\label{42:09} }
\Statex
\State{$\vect{w}_{k-1} = \vect{u}_{k-1} - {V}_{k-2}\:\vect{r}^{(2)}_{1:k-2,k-1}$ 
\Comment{Jacobi iteration} \label{42:10} }
\State{$\vect{v}_{k-1} = 
\: \vect{w}_{k-1} \:/ \: {\gamma_{k-1}}$ \label{42:11} }
\Statex
\State{$\vect{r}^{(1)}_{1:k-1,k} =
\left[ \begin{array}{cc} I & 0  \\
-L_{k-1,1:k-2} & 1 \end{array} \right]
\left[ \begin{array}{c} \vect{r}^{(0)}_{1:k-2,k} \\  
\vect{r}^{(0)}_{k-1,k}\end{array} \right]$ }
\State{$\vect{u}_k  = A\vect{u}_{k-1} \: / \: \gamma_{k-1} - 
{V}_{k-1}\:\vect{r}^{(1)}_{1:k-1,k}$ \Comment{Gauss-Seidel iteration} \label{42:13} }
\Statex
\State{$H_{1:k-2,k-1} =  
\vect{r}^{(4)}_{1:k-2,k-1} +  \vect{r}^{(3)}_{1:k-2,k-1}$
\Comment{Arnoldi expansion $H_{k-1}$} \label{42:14} }
\State{$\vect{r}^{(4)}_{1:k-1,k} = \vect{r}^{(0)}_{1:k-1,k}
- \ H_{1:k-1,1:k-2}\:\vect{r}^{(3)}_{1:k-2,k-1}$  \ \label{42:15} }
 \Statex
\EndFor
\State{$\vect{y}_m = {\rm argmin}  \:\|\: H_{m}\vect{y}_m - \rho\:\vect{e}_1 \:\|_2$ ,
 $\vect{x}_m = \vect{x}_0 + V_m\vect{y}_m$}
\end{algorithmic}
\caption{\label{onereduce} Low-Synch One-Reduce hybrid MGS-CGS GMRES}
\end{algorithm}

The triangular solve and matrix-vector multiply 
for the multiplicative iterations require $(k-1)^2$ flops at
iteration $k$ and thus lead to a slightly higher operation count 
versus the traditional MGS algorithm, which is $2m^2n$ for an $n\times m$ matrix.   
The  matrix-vector multiplies in Step \ref{42:04} of Algorithm \ref{onereduce}  
have complexity $4nk$ at iteration $k$ and the norm in Step \ref{42:06}
requires $2n$ flops, for a total of $4mn^2 + 3/3 n^3 + {\cal O}(mn)$ flops.
The number of global reductions is decreased from $k-1$ at iteration $k$ in MGS-GMRES
to only one when combined with the lagged normalization of the Krylov basis 
vectors. These costs are comparable to the DCGS-2 algorithm 
requiring $4mn^2$ flops. Indeed, the MGS-CGS GMRES
can be viewed as a hybrid algorithm  with Gauss-Seidel (MGS) and
Jacobi (CGS) iterations.

The traditional MGS-GMRES algorithm computes an increasing number
of inner products at each Arnoldi iteration. These take the
form of global reductions implemented as {\tt MPI\_Allreduce}.
A global reduction requires ${\cal O}(\log P)$ time to complete,
where $P$ is the number of MPI ranks running on a parallel machine.
This can be further complicated by heterogeneous computer architectures
based on graphical processor units (GPUs). However, the
single GPU performance of DCGS-2 is well over 200 GigaFlops/sec,
and merges the matrix-vector products in Steps \ref{42:10} and \ref{42:13} above
into one GPU kernel call for increased execution speed.
Recent improvements to MPI collective global communication
operations, that reduce latencies, include node-aware optimizations as
described by Bienz et al.~\cite{BienzGroppOlson2019}.
Among the different parallel variants of the algorithms
studied by Yamazaki et al.~\cite{2020-yamazaki-proceedings-of-siam-pp20},
the low-synch implementation of MGS-GMRES exhibited the best strong scaling 
performance on the ORNL Summit supercomputer.

\subsection*{Backward stability}\label{sec:depNorm}

The MGS-GMRES algorithm was proven to be backward stable for the solution 
of linear systems $A\vect{x} = \vect{b}$ 
in \cite{2006--simax--paige-rozloznik-strakos} and 
orthogonality is maintained to ${\cal O}(\eps)\kappa({B}_k)$, where
$B_k = [\vect{r}_0,\:A \bar{V}_k]$
and $\tv_k$ is the matrix  generated in finite precision arithmetic 
at iteration $k$,
with computed vectors as columns.
Our backward error analysis of the iterated Gauss-Seidel Gram-Schmidt 
Algorithm \ref{fig:PSGSA} is also based on the  $QR$ factorization of the matrix $\bar{B}_k$. 
The error matrix, $F_k$, for the computed Arnoldi expansion 
after $k$ iterations is expressed as 
\[
A\tv_k - \tv_{k+1}\bar{H}_{k+1,k} = F_k, 
\]
and is a matrix that grows in size by one column at each iteration.
  Recall that 
the strictly lower triangular matrix $L_{k}$
is incrementally computed  one row
per iteration as in \eqref{eq:katvec} and is obtained from the relation
$
\tv_{k}^T\tv_{k} = I + L_{k} + L_{k}^T.
$
Thus it follows from the analysis in Section 3 that the loss 
of orthogonality is 
\begin{equation}
\|\:I - \tv_{k}^T\tv_{k}\:\|_2 \lesssim {\cal O}(\eps). 
\end{equation}
The IGS-GMRES algorithm (Algorithm \ref{lowsynch}) thus 
maintains orthogonality to the level ${\cal O}(\eps)$. 
Therefore, it follows from Drko\v{s}ov\'{a} et al.~\cite{drsb1995},
that under the reasonable assumption of the numerical non-singularity 
of the coefficient matrix, the algorithm is also backward stable. 

\begin{figure}
\centering
\includegraphics[width=0.7\textwidth]{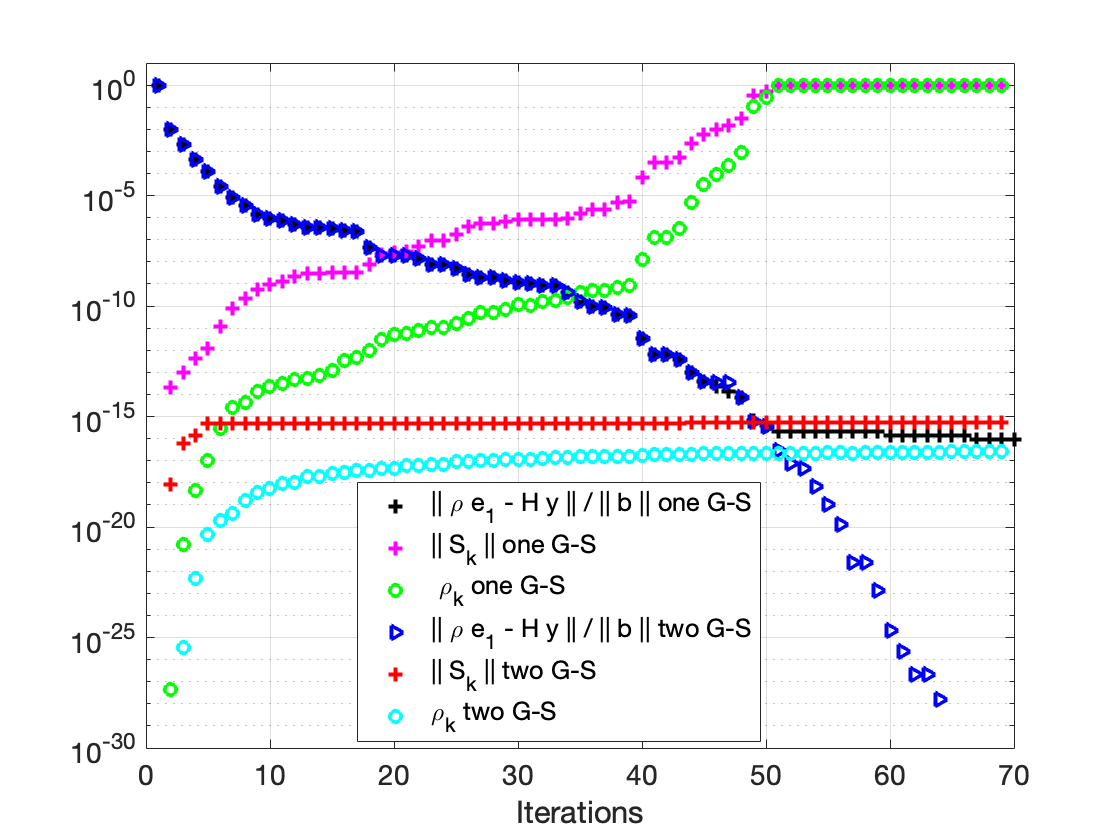}
\caption{\label{fig:f1836-bounds}fs1836 matrix. Arnoldi relative residual, spectral radius of $M_k^{-1}N_k$ and Paige's $\|S_k\|$ for one and two Gauss-Seidel iterations.}
\end{figure}

Our first experiment illustrates that the bounds derived in the previous 
sections properly capture the behavior of the IGS-GMRES algorithm in finite precision. 
In the next section, we demonstrate that our bounds 
continue to hold for a broad class of matrices.
In particular, we examine the {\tt fs1836} matrix studied by Paige and Strako\v{s} 
\cite{Paige2002}. Our Figure \ref{fig:f1836-bounds} should be compared with
Figures 7.1 and 7.3 of their (2002) paper. In order to demonstrate empirically 
that the backward error is reduced by the iteration matrix
$M_{k-1}^{-1}N_{k-1}$, the quantity $\|S_k^{(2)}\|_2$ is computed,
as defined by Paige et al.~\cite{2006--simax--paige-rozloznik-strakos},
which  measures the loss of orthogonality for two iterations.
The spectral radius $\rho_k$ of the matrix $M_k^{-1}N_k$
is highly correlated with, and follows the metric, $\|S_k^{(2)}\|_2$.
The metric $\|S_k^{(2)}\|_2$ is plotted for one and two Gauss-Seidel
iterations in Figure \ref{fig:f1836-bounds}. The Arnoldi 
relative residuals are plotted as in Figure 7.1 of 
Paige and Strako\v{s} \cite{Paige2002}, which, for one iteration,
stagnate near iteration forty-three at $1\times 10^{-7}$ 
before reaching ${\cal O}(\eps)$. For two Gauss-Seidel iterations,
the (Arnoldi) relative residual continues to decrease monotonically
and the norm-wise relative backward error (\ref{eq:nrbe}) reaches the level of machine precision: $\beta(\vect{x}_k) = 6.6\times 10^{-17}$ at iteration fifty.
Most notably, the 2-norm of $A$ is large, where $\|A\|_2 = 1.2\times 10^9$.

A normal matrix $X \in \mathbb{C}^{n\times n}$ satisfies 
$X^*X = XX^*$. 
Henrici's definition 
of the departure from normality 
\begin{equation}
    {\rm dep}(X) = \sqrt{\|X\|_F^2-\|\Lambda(X)\|_F^2},
\end{equation}
where $\Lambda(X) \in\mathbb{C}^{n \times n}$ is the diagonal matrix 
containing the eigenvalues of $X$ \cite{Henrici1962}, serves 
as a useful metric for the loss of orthogonality.  
While we find practical use for this metric for
measuring the degree of (non)normality of a matrix, 
there are of course other useful metrics to 
describe (non)normality. We refer the reader to
\cite{Ipsen1998,Henrici1962,Trefethen2005} 
and references therein. 
In particular, the loss of orthogonality 
is signaled by the departure from  normality as follows:
\begin{equation}
{\rm dep}^2(M_k^{-1}) =  {\rm dep}^2(I + L_k) = 
\|I + L_k\|_F^2 - \|I\|_F^2 = \|I\|_F^2 + \|L_k\|_F^2 - k = \|L_k\|_F^2.
\label{eq:LOO}
\end{equation}
Note that for the iteration matrix  $M^{-1}_k N_k= (I + L_k)^{-1} L_k^T$  
we have, and it is observed in practice, that 
$\|M^{-1}_k N_k\|_2 \leq  \|M^{-1}_k N_k\|_F \approx \|L_k\|_F$ up to the first order in $\eps$.  

\section{Numerical Results}

Numerically challenging test problems for GMRES have been proposed
and analyzed over the past 35 years. These include both
symmetric and non-symmetric matrices. 
Simoncini and Szyld \cite{2003_simoncini_SISC} 
introduced a symmetric, diagonal 
matrix with real eigenvalues, causing MGS-GMRES residuals to stagnate. 
Highly non-normal matrices from Walker \cite{Walker88} 
were used to explore the convergence characteristics of 
HH-GMRES and then the non-normal fs1836 from Paige et al. 
\cite{2006--simax--paige-rozloznik-strakos}
and west0132 from Paige and Strako\v{s}\cite{Paige2002}
encountering stagnation. In addition to these,
we consider 
the impcol\_e matrix from Greenbaum et al.~\cite{Green92}. Matrices with
complex eigenvalues forming a disc inside the unit circle such as the Helmert
matrix from Liesen and Tich\'y \cite{Liesen2004}, are also evaluated. 
Results from a  pressure continuity solver with AMG 
preconditioner and a circuit simulation with
the ADD32 matrix from Rozlo{\v{z}}n{\'{\i}}k, 
Strako\v{s}, and T\accent23uma 
\cite{RozloznikST96} are also presented. 
Numerical results were obtained with Algorithm 4.1
and we verified that Algorithm 4.2 computations are 
comparable. In particular, the diagonal
elements $R_{k,k}$ computed with the lagged norm
agree closely with those from Algorithm 4.1.
In all cases, the orthogonality between the computed vectors remains
on the level ${\cal O}(\eps)$, the 
smallest singular value of $\bar{V}_k$
remains close to one, and the Krylov basis vectors remain
linearly independent.
Therefore 
the norm of the true relative residual reaches the level ${\cal O}(\eps)$   at least 
by the final iteration, while
the Arnoldi residual continues to decrease far below this level.

\subsection*{Ill-Conditioned Diagonal Matrix}

\begin{figure}
\centering
\includegraphics[width=0.7\textwidth]{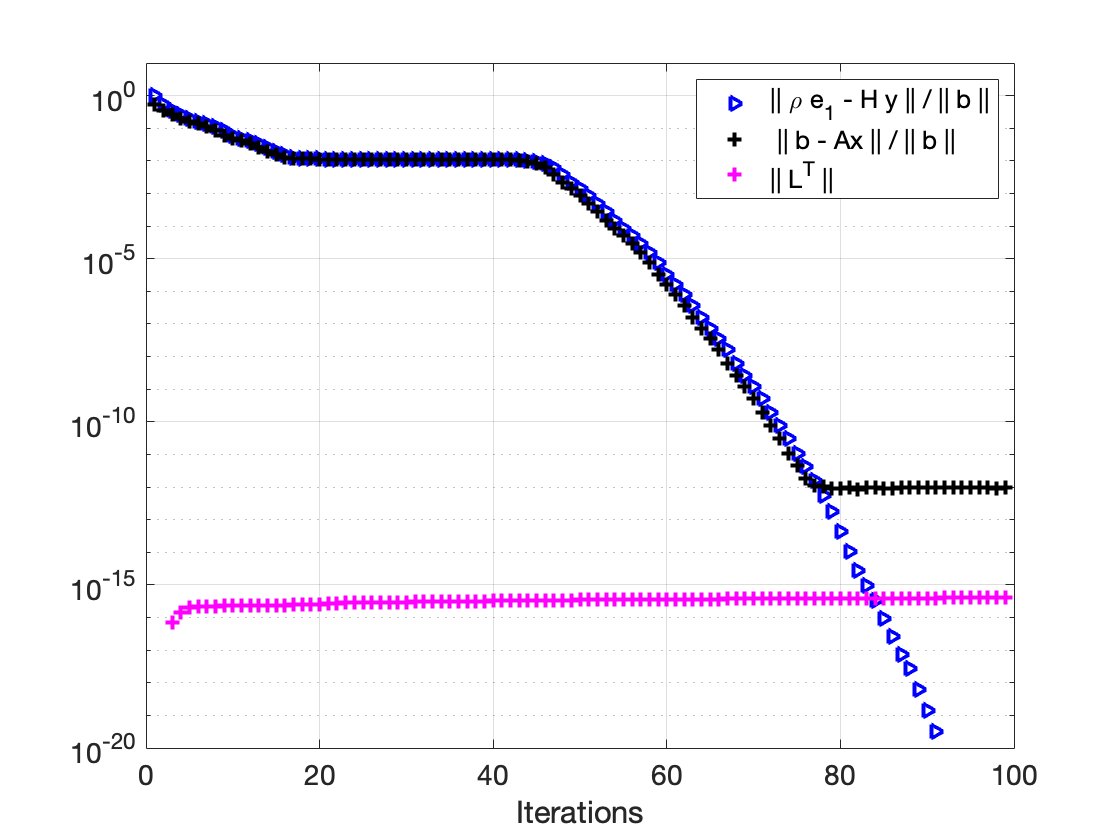}
\caption{\label{fig:simoncini} 
Simoncini matrix. Arnoldi relative residual.
Loss of orthogonality relation (\ref{eq:LOO}).}
\end{figure}

Simoncini and Szyld \cite{2003_simoncini_SISC} consider several 
difficult and very ill-conditioned problems that can lead to stagnation
of the residual before converging to 
the level of machine precision ${\cal O}(\eps)$.
In their example 5.5, they construct 
$A = {\rm diag}([1e-4, 2, 3, \ldots, 100])$, a diagonal matrix,
and the right-hand side $\vect{b} = {\rm randn}(100, 1)$ is 
normalized so that $\|\vect{b}\|_2 = 1$. The condition number
of this matrix is $\kappa(A)=1\times 10^6$ and $\|A\|_2 =100$.
With the MGS-GMRES algorithm, the relative residual stagnates at
the level $1\times 10^{-12}$ after 75 iterations, when $\|S_k\|_2 = 1$, 
indicating that the Krylov vectors are not linearly independent.
In the case of the iterative Gauss-Seidel formulation of GMRES 
(with two Gauss-Seidel iterations), the convergence history 
is plotted in Figure \ref{fig:simoncini}, where it can 
be observed that the relative Arnoldi residual continues to decrease
monotonically. Furthermore, the true relative residual
is plotted along with $\| L_{k-1} \|_F$. The latter indicates
that a significant loss of orthogonality does not occur.

\subsection*{Ill-Conditioned Symmetric and Non-Symmetric Matrices}

Figures 1.1 and 1.2 from Greenbaum et al.~\cite{Green92} 
describe the results for STEAM1 
(using the HH and MGS implementations of GMRES, respectively). Similarly, Figures 1.3 and 1.4
from Greenbaum et al.~\cite{Green92} correspond to IMPCOLE. 
They emphasize that the convergence behavior illustrated in
these plots is typical of the MGS-GMRES and HH-GMRES algorithms.
The condition number of the system matrix is 
$\kappa(A) = 2.855\times10^7$ and $\|A\|_2 = 2.2\times 10^7$ for STEAM1, whereas
$\kappa(A) = 7.102\times  10^6$ and $\|A\|_2 = 5.0\times 10^3$ for IMPCOLE.

\begin{figure}
\centering
\includegraphics[width=0.7\textwidth]{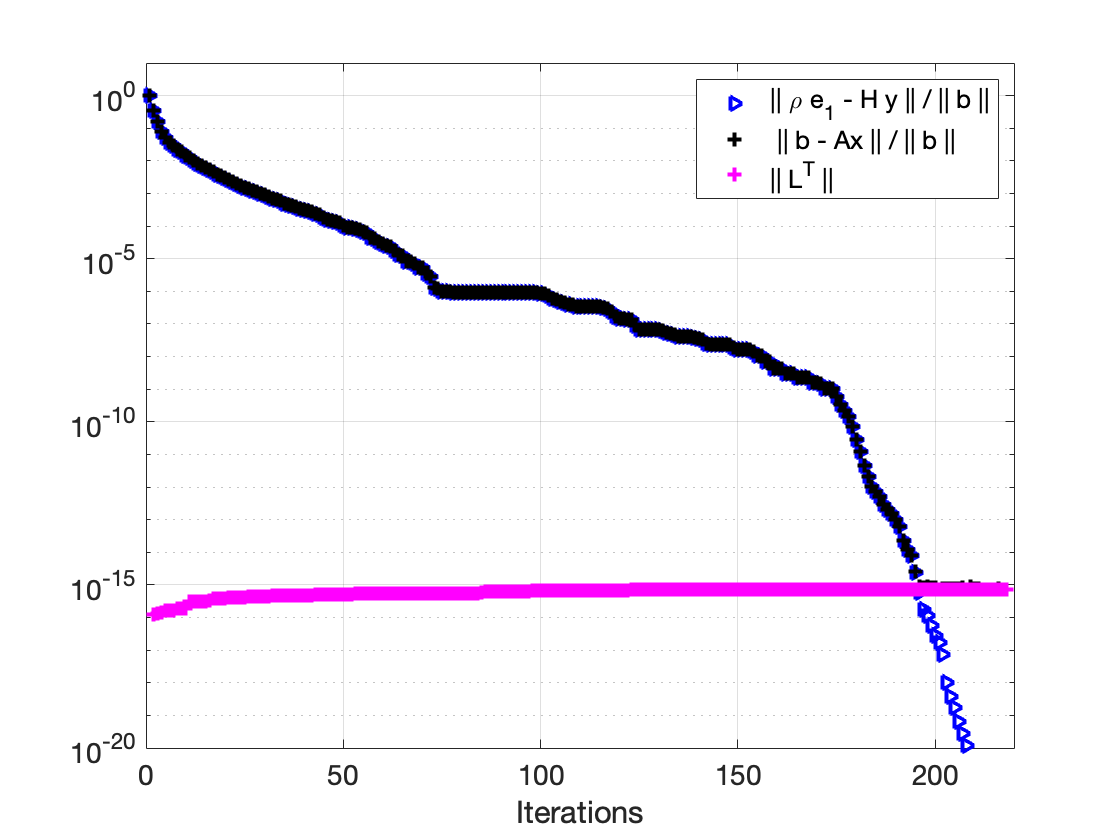}
\caption{\label{fig:steam1}
steam1 matrix. Arnoldi relative residual. 
Loss of orthogonality relation (\ref{eq:LOO}).}
\end{figure}

\begin{figure}
\centering
\includegraphics[width=0.7\textwidth]{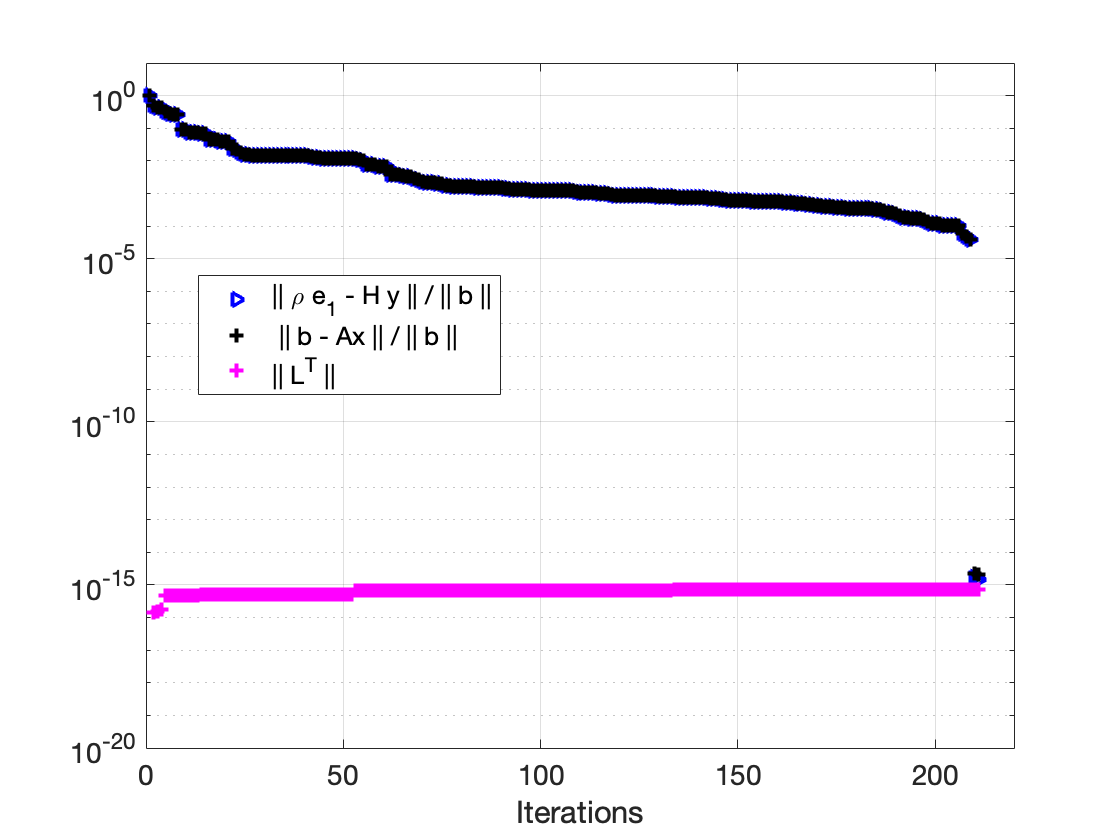}
\caption{\label{fig:impcolColScale}impcol\_e matrix.
Arnoldi relative residual. 
Loss of orthogonality relation (\ref{eq:LOO}).}
\end{figure}

Greenbaum et al.~\cite{Green92} observe that 
although orthogonality of the Krylov vectors  
is not maintained near machine precision, 
as is the case for the Householder implementation,
the relative residuals of the MGS-GMRES algorithm are almost identical 
to those of the HH-GMRES until the smallest singular value of the matrix $\bar{V}_k$
begins to depart from one. At that point the MGS-GMRES relative residual norm
begins to stagnate close to its final precision level. This observation is 
demonstrated with the numerical examples for matrices STEAM1 ($N = 240$, 
symmetric positive definite matrix used in oil recovery simulations) and 
IMPCOLE ($N = 225$, nonsymmetric matrix from modelling of the hydrocarbon 
separation problem).
In both experiments $\vect{x} = (1,..., 1)^T$, 
$\vect{b} = A\vect{x}$ and $\vect{x}_0 = {\bf 0}$. 
The convergence histories for the iterated Gauss-Seidel GMRES algorithm applied to
these matrices are plotted in Figures \ref{fig:steam1} and 
\ref{fig:impcolColScale}. A significant loss of orthogonality 
is not observed until the last iteration.
Otherwise the computed metric $\|L_{k-1}\|_F$ and the true
relative residual
remain near ${\cal O}(\eps)$. The Arnoldi residual continues
to decrease and
the smallest singular value of $V_k$ is $0.99985$.

\subsection*{Highly Non-Normal Matrices}

\begin{figure}
\centering
\includegraphics[width=0.7\textwidth]{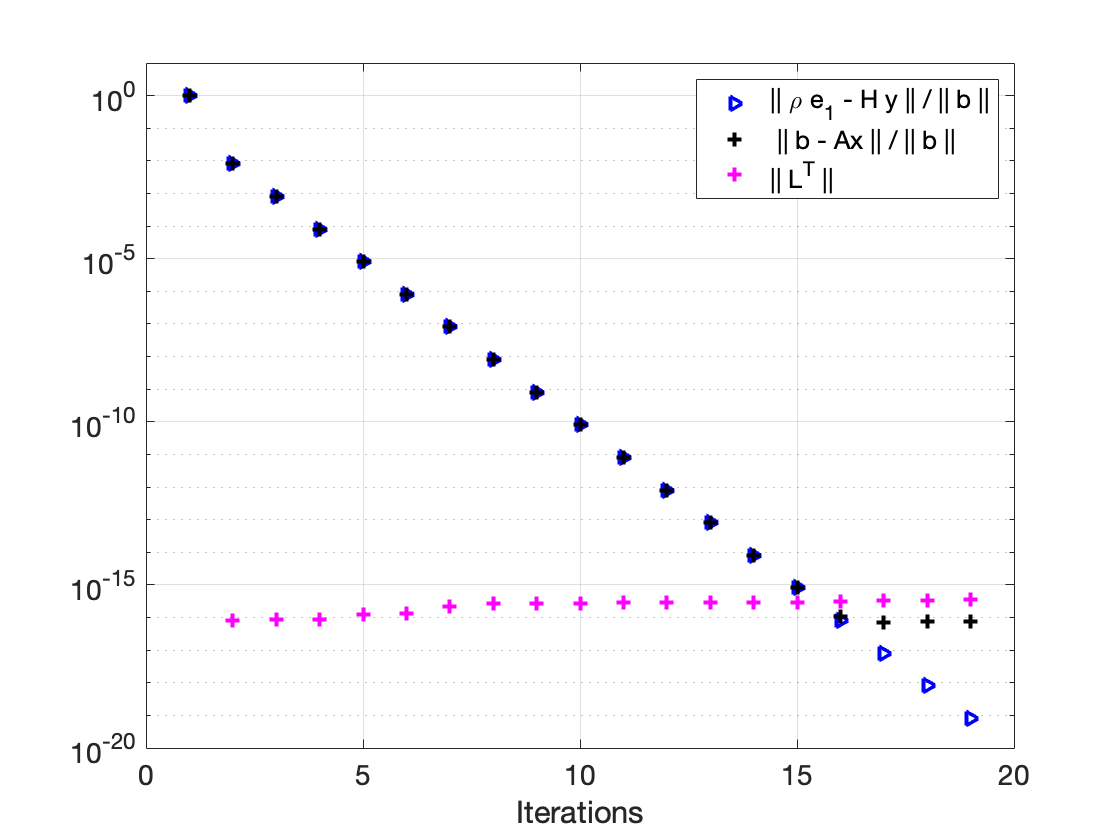}
\caption{\label{fig:Embree} 
Embree $\delta$  matrix. Arnoldi relative residual.
Loss of orthogonality relation (\ref{eq:LOO}).}
\end{figure}

Bidiagonal matrices with a $\delta$ off-diagonal were
studied by Embree \cite{Embree99}. These are non-normal
matrices where $0 < \delta \le 1$ and are completely 
defective for all $\delta \ne 0$. A defective
matrix is a square matrix 
that does not have a complete basis of eigenvectors, 
and is therefore not diagonalizable, and
the pseudo-spectra \cite{Trefethen2005} of these
matrices are discs in the complex plane. Our
IGS-GMRES algorithm leads to convergence after 16 iterations
without stagnation and orthogonality is maintained
to machine precision as plotted in Figure \ref{fig:Embree}.
The matrix 2-norm is $\|A\|_2 = 1.1$ and condition $\kappa(A) = 1.2$.
Walker \cite{Walker88} employed the highly non-normal matrix in
equation (\ref{eq:walker})  to compare the Gram-Schmidt 
and Householder implementations
of GMRES. The element $\alpha$ controls both the condition number $\kappa(A)$
and the departure from normality ${\rm dep}(A)$ of the matrix of size $n\times n$. 
Here $\|A\|_2 = 2.0\times 10^3$.
\begin{equation}
A = \left[
\begin{array}{ccccc}
 1  &  0 & \cdots & 0 & \alpha\\
 0  &  2 & \cdots & 0 & 0 \\
\vdots & \vdots &  & \vdots & \vdots \\
0 & 0 & \cdots & 0 & n
\end{array}
\right], \quad
\vect{b} = \left[
\begin{array}{c}
1 \\
1 \\
\vdots \\
1
\end{array}
\right]
\label{eq:walker}
\end{equation}
For large values of $\alpha$, Walker found that the MGS
residual would stagnate and that the CGS algorithm led to instability.
Furthermore, it was found that even CGS-2 with re-orthogonalization
exhibited some instability near convergence. 
HH-GMRES maintains ${\cal O}(\eps)$ orthogonality as measured 
by $\|I - \bar{V}_k^T\bar{V}_k\|_F$ and reduces the 
relative residual to machine precision.

\begin{figure}
\centering
\includegraphics[width=0.7\textwidth]{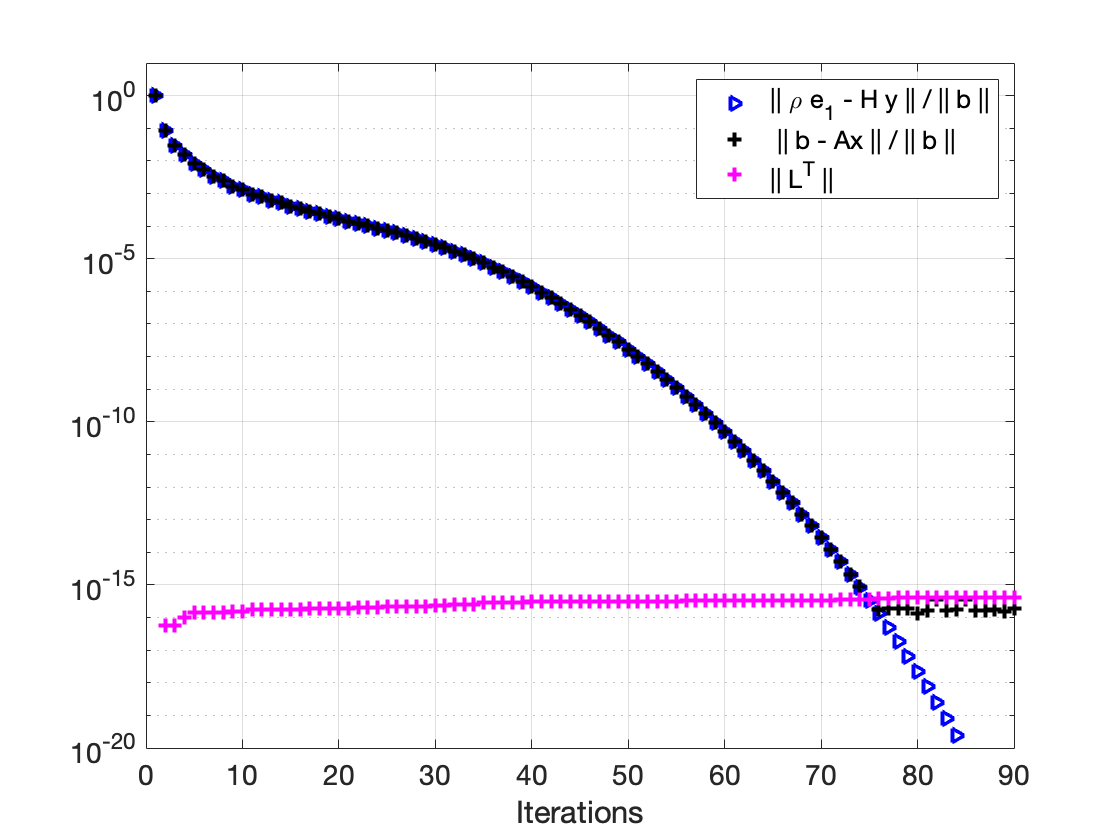}
\caption{\label{fig:WalkerrColScale}Walker matrix.
Arnoldi relative residual.
Loss of orthogonality relation (\ref{eq:LOO}).}
\end{figure}

In our experiments, the value $\alpha = 2000$ leads to a
matrix with $\kappa(A) = 4\times 10^5$. The departure from normality,
based on Henrici's metric, is large: ${\rm dep}(A) = 2000$.
The convergence history is displayed in
Figure \ref{fig:WalkerrColScale}. The loss of orthogonality
remains near $O(\eps)$ and our upper bound is close for this problem.
Notably, the oscillations present in the relative
residual computed by the CGS-2 variant
are not present in the iterated Gauss-Seidel convergence history plots, 
where the Arnoldi relative residual decreases monotonically.

\begin{figure}
\centering
\includegraphics[width=0.7\textwidth]{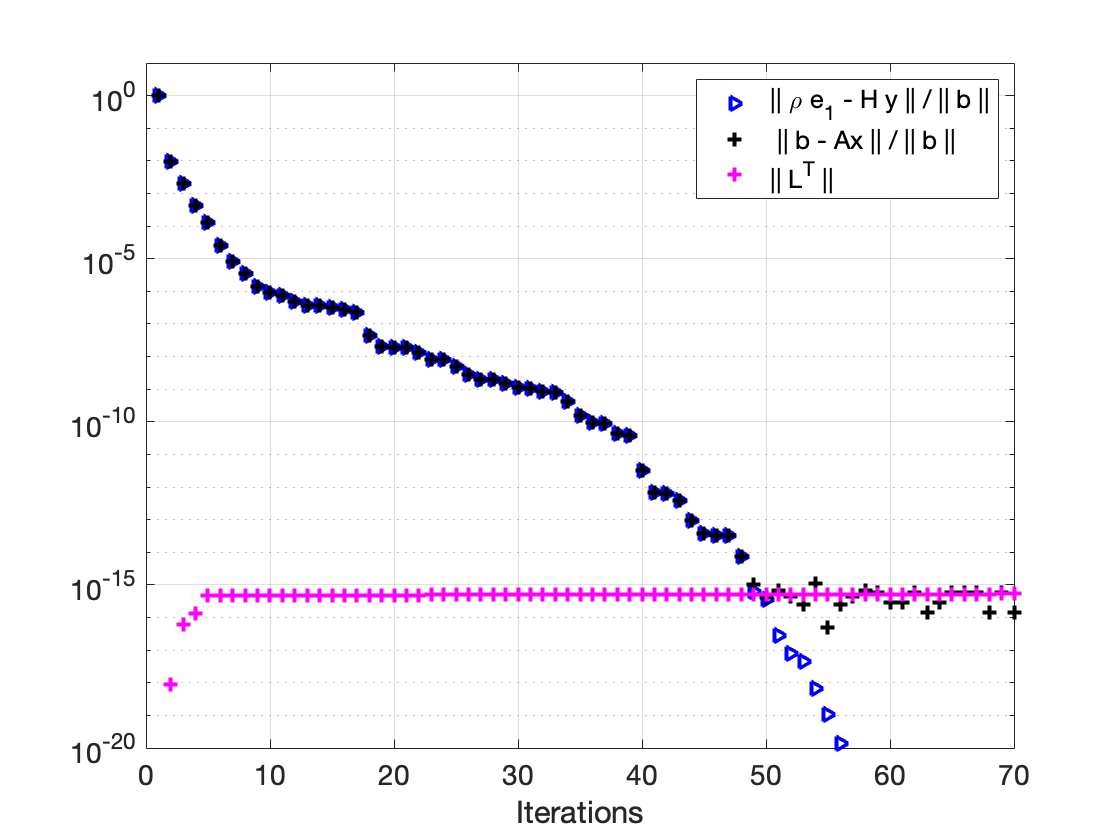}
\caption{\label{fig:f1863ColScale}
fs1863 matrix. Arnoldi relative residual. 
Loss of orthogonality relation (\ref{eq:LOO}).}
\end{figure}

\begin{figure}
\centering
\includegraphics[width=0.7\textwidth]{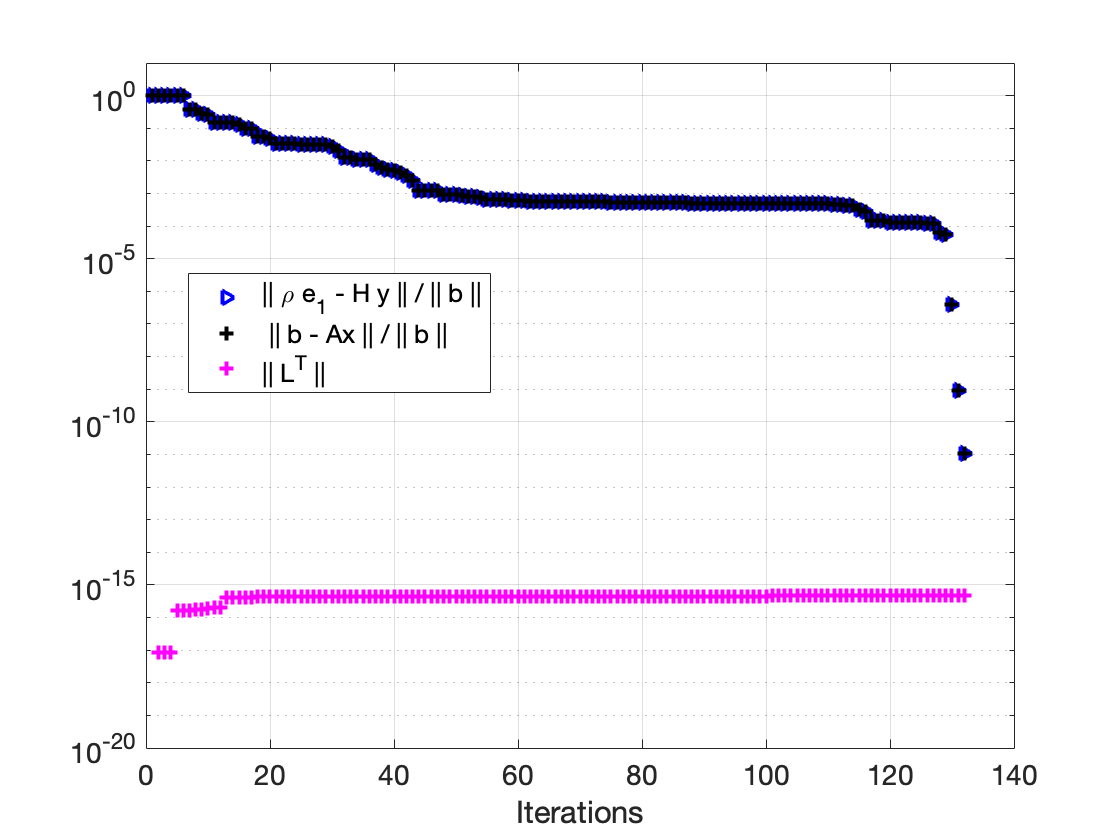}
\caption{\label{fig:westColScale}west0132 matrix.
Arnoldi relative residual. 
Loss of orthogonality relation (\ref{eq:LOO}).}
\end{figure}

Paige and Strako\v{s} \cite{Paige2002} experimented with MGS-GMRES  convergence 
for the non-normal matrices  
FS1836 and WEST0132. In all their experiments $\vect{b} = (1, \ldots, 1)^T$ .
The matrix FS1836 has dimension $n = 183$, with   $\|A\|_2 \approx 1.2\times 10^9$, and 
$\kappa(A) \approx 1.5\times 10^{11}$. The matrix WEST0132 has $n = 132$, 
$\|A\|_2 \approx 3.2\times 10^5$, and $\kappa(A) \approx 6.4\times 10^{11}$.
Their Figure 7.1 indicates that the relative residual for FS1836
stagnates at $1\times 10^{-7}$ at iteration 43 when orthogonality is lost.
The relative residual for the WEST0132 matrix also stagnates at
the $1\times 10^{-7}$ level after 130 iterations. These results contrast with 
our Figures \ref{fig:f1863ColScale}
and \ref{fig:westColScale}. In both cases the Arnoldi residuals 
continue to decrease and $\|L_{k-1}\|_F$ grows
slowly or remains close to machine precision. 
The smallest singular value remains at $\sigma_{\min}(\bar{V}_k) = 1$
for both matrices.

\subsection*{Complex Eigenvalues in a Disc}

\begin{figure}
\centering
\includegraphics[width=0.7\textwidth]{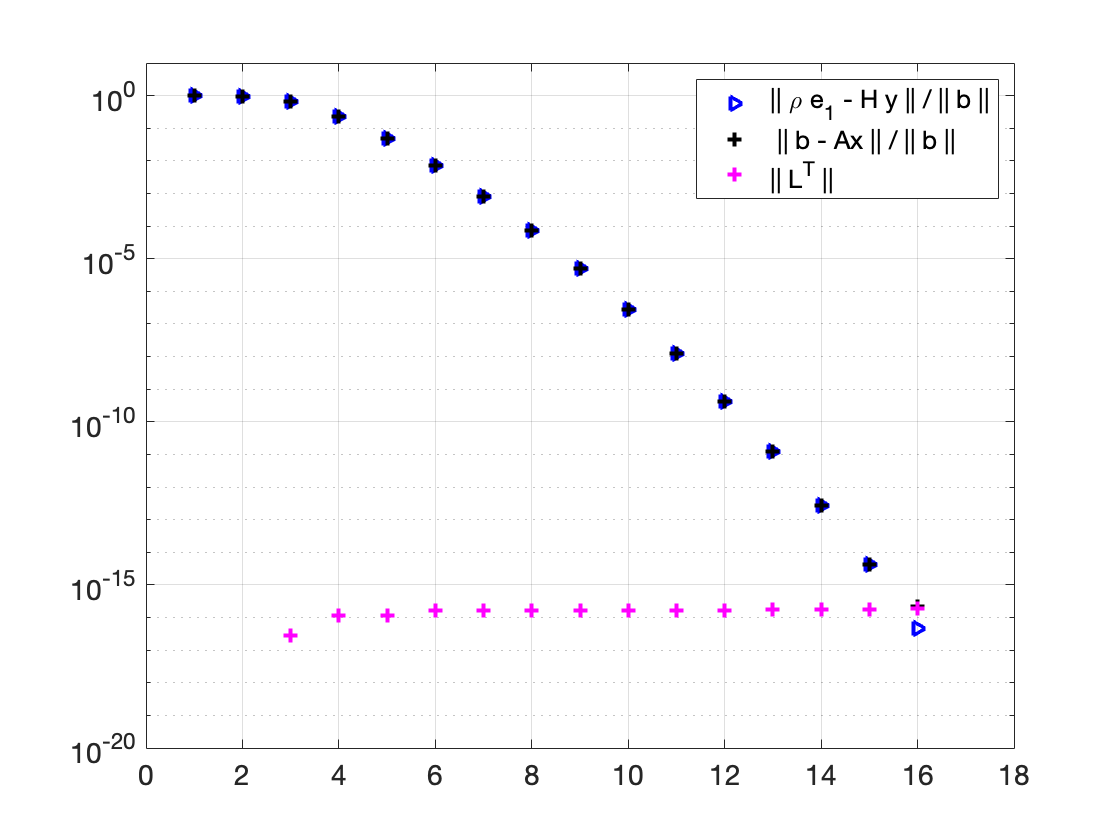}
\caption{\label{fig:Helmert2}
Helmert matrix. Arnoldi relative residual. 
Loss of orthogonality relation (\ref{eq:LOO}).}
\end{figure}

Liesen and Tich\'{y} \cite{Liesen2004} employ
the Helmert matrix generated by the MATLAB  
command {\tt gallery('orthog',18,4)}. Helmert matrices occur in a number of 
practical problems, for example in applied statistics. The matrix is 
orthogonal, and the eigenvalues cluster around $-1$, as in the right panels of 
their Figure 4.4. The worst-case MGS-GMRES residual norm decreases quickly throughout 
the iterations and stagnates at the 12--th iteration, where the 
relative residual remains at $1\times 10^{-10}$.
From the convergence history plotted in Figure \ref{fig:Helmert2},
the loss of orthogonality remains near 
machine precision and the Arnoldi relative residual does not stagnate. The
quantity $\|L_{k-1}\|_F$
is an excellent predictor of the orthogonality. 

\subsection*{Nalu-Wind Model}

Nalu-Wind solves the incompressible Navier-Stokes
equations, with a pressure projection.  
The governing equations are discretized in time with a BDF-2 integrator,
where an outer Picard fixed-point iteration is employed to reduce the nonlinear
system residual at each time step.
Within each time step, the Nalu-Wind simulation time is often dominated by the
time required to setup and solve the linearized governing equations.
The pressure systems are solved using GMRES with an
AMG preconditioner, where a polynomial Gauss-Seidel
smoother is now applied are described in Mullowney 
et al.~\cite{Mullowney2021}.  Hence, relaxation is a compute time intensive
component, when employed as a smoother.

The McAlister experiment for wind-turbine blades is an
unsteady RANS simulation of a fixed-wing, with a NACA0015 cross section, 
operating in uniform inflow. 
Resolving the high-Reynolds number boundary layer over the wing surface 
requires resolutions of ${\cal O}(10^{-5})$ normal to the surface resulting in 
grid cell aspect ratios of ${\cal O}(40,000)$. These high aspect ratios present a 
significant challenge.
The simulations were performed for the wing at 12 degree angle of
attack, 1 m chord length, 
denoted $c$, aspect ratio of 3.3, i.e., $s = 3.3c$, and square wing tip. 
The inflow velocity is $u_{\infty}= 46$ m/s, the density is $\rho_{\infty} = 1.225$
${\rm kg/m}^3$, and dynamic viscosity is $\mu = 3.756\times 10^{-5}$ kg/(m s), 
leading to a Reynolds number $Re = 1.5\times 10^6$.
Due to the complexity of mesh generation, 
only one mesh with approximately 3 million grid points was generated. 

\begin{figure}
\centering
\includegraphics[width=0.7\textwidth]{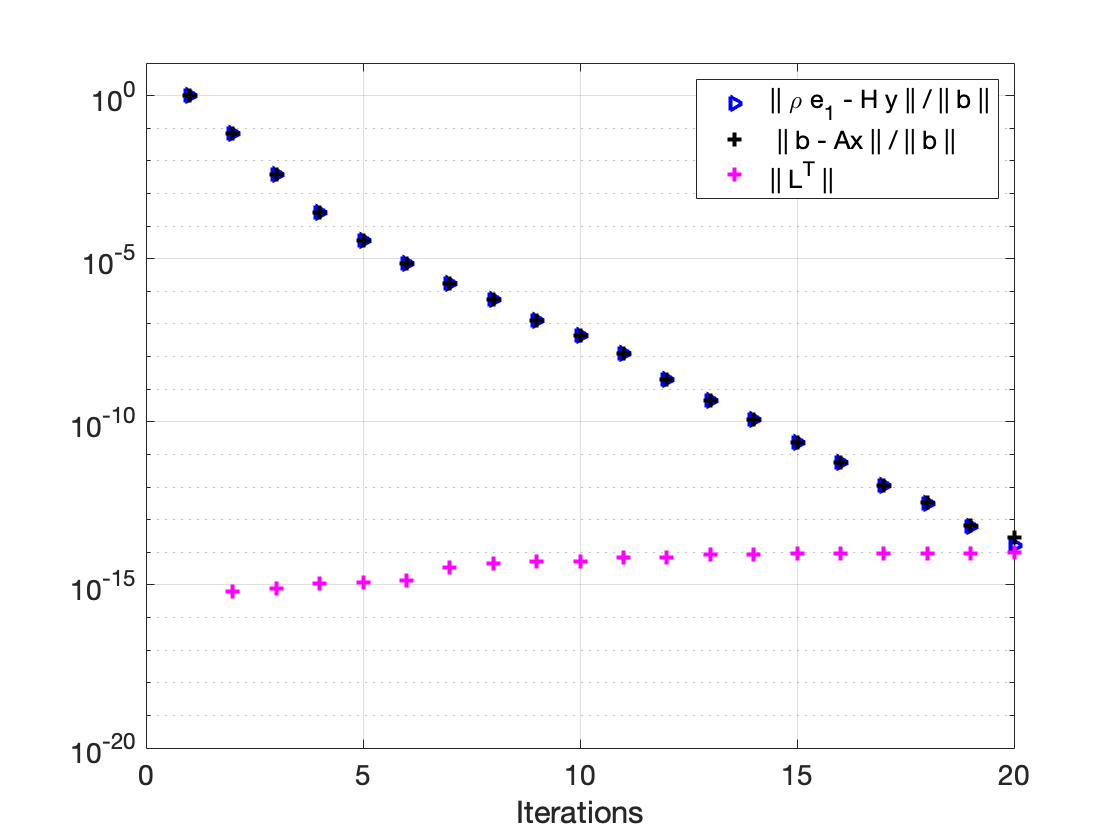}
\caption{\label{fig:Pressure}
Pressure matrix. Arnoldi relative residual. 
Loss of orthogonality relation (\ref{eq:LOO}).}
\end{figure}

The smoother is hybrid block-Jacobi with two sweeps of polynomial
Gauss-Seidel relaxation applied locally on a subdomain and then Jacobi 
smoothing for globally shared
degrees of freedom.  The coarsening rate for the wing simulation is roughly
$4\times$ with eight levels in the $V$-cycle for {\it hypre}
\cite{Falgout2004}. Operator complexity
$C$ is approximately $1.6$ indicating more efficient $V$-cycles with aggressive
coarsening, however, an increased number of solver iterations 
are required compared to standard coarsening.  The convergence
history is plotted in Figure \ref{fig:Pressure}, where the 
loss of orthogonality is completely flat and close to machine precision.

\subsection*{Circuit Simulation}

\begin{figure}
\centering
\includegraphics[width=0.7\textwidth]{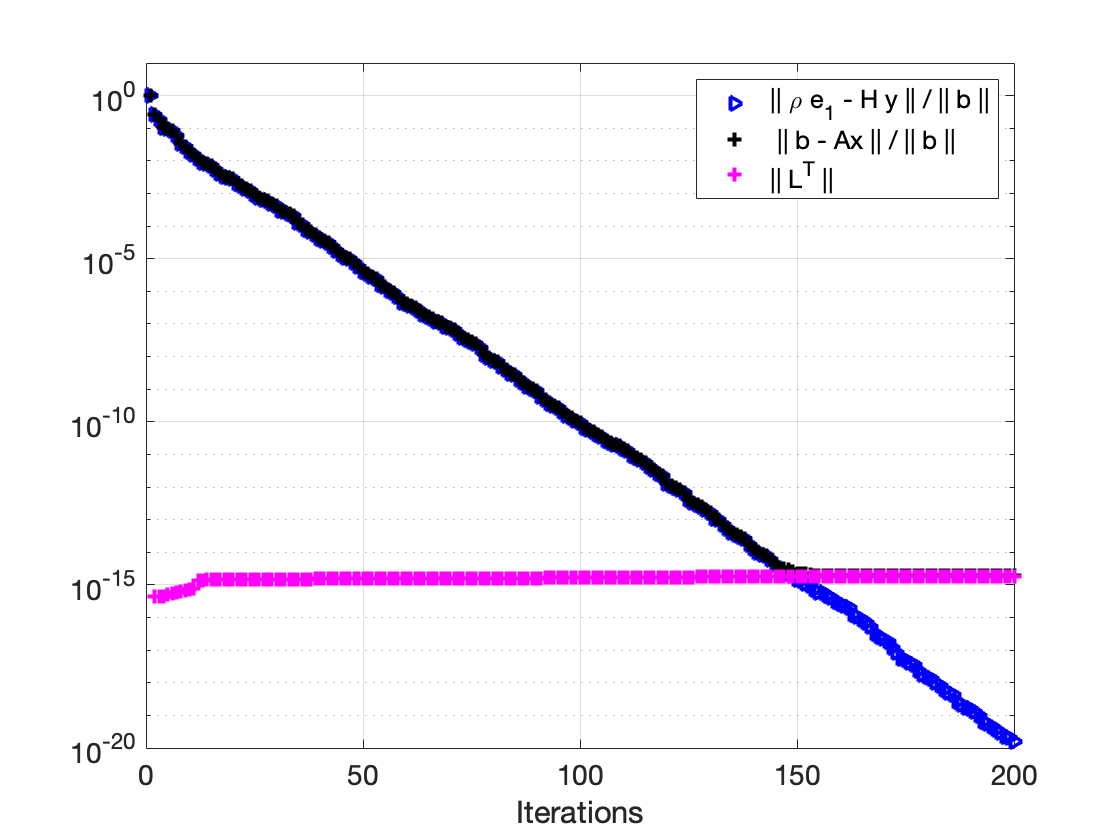}
\caption{\label{fig:Add32}Add32 matrix. Arnoldi relative residual. 
Loss of orthogonality relation (\ref{eq:LOO}).}
\end{figure}

Rozlo{\v{z}}n{\'{\i}}k et al.~\cite{RozloznikST96} study a typical linear
system arising in circuit simulation (the matrix from a 32-bit adder design). The matrix has $\|A\|_2 = 0.05$ and $\kappa(A) = 213$. 
In exact arithmetic the Arnoldi vectors are orthogonal. However, in finite
precision computation the orthogonality is lost, which may potentially affect 
approximate solution. In their Figure 3, the authors have plotted the loss of
orthogonality of the computed Krylov vectors for different implementations of 
the GMRES method (MGS, Householder and CGS). The comparable results for
IGS-GMRES are plotted in Figure \ref{fig:Add32}. 
The smallest singular value remains at $\sigma_{\min}(\bar{V}_k) = 1$.

\section{Conclusions}

The essential contribution of our work was to derive an iterative Gauss-Seidel
formulation of the GMRES algorithm due 
to Saad and Schultz \cite{Saad86}
that employs approximate solutions of the normal
equations appearing in the Gram-Schmidt projector, 
based on observations of Ruhe \cite{Ruhe83} and the low-synchronization 
algorithms introduced by \'{S}wirydowicz et al.~ \cite{2020-swirydowicz-nlawa}.


The insights gained from the seminal work of Ruhe \cite{Ruhe83} led us to 
the conclusion that the modified Gram-Schmidt algorithm is equivalent to
one step of a {\it multiplicative} Gauss-Seidel iteration method applied to the normal equations
$Q_{k-1}^TQ_{k-1}\:\vect{r}_{1:k-1,k}= Q_{k-1}^T\vect{a}_k$.
Similarly, the classical Gram-Schmidt algorithm can be viewed as one step of an {\it additive}
Jacobi relaxation. The projector is then given by 
$P\vect{a}_k = \vect{a}_k - Q_{k-1}\:T_{k-1}\:Q_{k-1}^T\vect{a}_k$,
where $T_{k-1}$ is a correction matrix. In the case of DCGS-2,
with delayed re-orthogonalization, Bielich et al.~\cite{DCGS2Arnoldi}
split and apply the symmetric (normal) correction matrix across two 
Arnoldi iterations and then apply Stephen's 
trick to maintain orthogonality.
For MGS, the lower triangular matrix $T_{k-1} \approx(\:Q_{k-1}^TQ_{k-1}\:)^{-1}$ 
appearing in \'{S}wirydowicz et al. \cite{2020-swirydowicz-nlawa} was identified as the 
inverse compact $WY$ form with $T^{(1)}_{k-1} = ( I + L_{k-1} )^{-1}$, where the
strictly lower triangular matrix $L_{k-1}$ was computed from the loss of 
orthogonality relation 
\[
Q_{k-1}^T Q_{k-1} = I + L_{k-1} + L_{k-1}^T.
\]
The matrix $T_{k-1}$ from the inverse compact $WY$ form of
the Gram-Schmidt projector was also present, 
without having been explicitly defined,
in the rounding error analysis of Bj\"orck \cite{Bjork_1967}, 
in Lemma 5.1.
In effect, the low-synchronization Gram-Schmidt 
algorithm presented in \cite{2020-swirydowicz-nlawa} 
represents one iteration step, with a zero initial guess, 
to construct an approximate projector. 
When two iteration steps are applied, the resulting correction
matrix $T_{k-1}^{(2)}$ is close to a symmetric (normal) matrix: 
\[
T_{k-1}^{(2)} = M_{k-1}^{-1}\: [ \: I + N_{k-1}\:M_{k-1}^{-1} \: ] =
T_{k-1}^{(1)} - T_{k-1}^{(1)}\:L_{k-1}^T\:T_{k-1}^{(1)}.
\]
However, the Gauss-Seidel formulation described by Ruhe \cite{Ruhe83}
differs from MGS-2 in floating-point arithmetic and the incremental MGS-2 
iterative refinement steps allow us to prove backward stability. 
When employed to compute the Arnoldi-$QR$ expansion, 
the GMRES formulation with two Gauss Seidel iterations
results in an ${\cal O}(\eps)$ backward error, preserving the orthogonality of the 
Krylov basis vectors to the level ${\cal O}(\eps)$,  
which is measured by $\|L_{k-1}\|_F$.
This result is related to recent work on the iterative solution 
of triangular linear systems using Jacobi iterations that
may diverge for highly non-normal triangular matrices 
\cite{CHOW2018219}. 

Here, the  departure from normality of the inverse of the 
correction matrix $T^{(1)}_{k-1}$ is a measure
of the loss of orthogonality. In this formulation, 
the matrix is lower triangular
and can substantially depart from normality as signaled by
$\sigma^2_{\max}(I + L_{k-1}) > 1$. Because the correction matrix $T^{(2)}_{k-1}$
associated with two iterations of Gauss--Seidel is close
to a symmetric (normal) matrix, the singular values of the
inverse remain close to one.
A departure from normality indicates a possible loss of 
numerical rank for the Krylov basis  
vectors $V_k$ with the smallest singular value decreasing from one. 
Our numerical
experiments, on challenging problems proposed over the past thirty-five 
years, demonstrate the robust convergence and low relative error
achievable with our IGS-GMRES.
Furthermore, with one iteration step the loss of
orthogonality is at most ${\cal O}(\eps)\kappa(B_k)$ 
and for the algorithm with two iterations, 
it remains near machine precision.

We have demonstrated that the iterated Gauss Seidel GMRES with two iterations
is backward stable and does not exhibit stagnation in the Arnoldi residual as defined
by Greenbaum et al.~\cite{Green92}. This refers to the approximate
residual norm computed as
$\| \rho\: \vect{e}_1 - H_{k+1,k}\vect{y}_k\:\|_2 / \rho$,
where $\vect{x}_k = \vect{x}_0 + V_k\vect{y}_k$ is the approximate solution
and the true residual is given by $\vect{b} - A\vect{x}_k$. The Arnoldi residual
does not stagnate and will continue to decrease monotonically. In particular,
the true relative residual will decrease towards machine 
precision level and as noted by Paige et al.~\cite{2006--simax--paige-rozloznik-strakos}, 
the norm-wise relative backward error $\beta(\vect{x}_k)$ will be ${\cal O}(\eps)$. 
In the original algorithm, the Arnoldi residual can stagnate when the
Krylov vectors lose linear independence and the smallest singular
value of $V_k$ decreases towards zero. For the iterative Gauss-Seidel GMRES algorithm
described herein, the singular values remain close to one
and $\beta(\vect{x}_k)$ remains close to 
machine precision at convergence. 

We have also presented a one-reduce algorithm which can be interpreted as an MGS-CGS  hybrid. 
This algorithm merits further study. Algorithm \ref{lowsynch} is equivalent to
the MGS-CGS GMRES when $L_{k-2} = 0$ at step $k-1$ 
in exact arithmetic and produces similar
results in finite precision (in terms of the loss of orthogonality and Arnoldi residual). 
Our numerical experiments demonstrated that the computed results
are comparable to the two-reduce IGS-GMRES
and diagonal elements $R_{k,k}$ agree closely
when the lagged norm is employed.
In particular, the loss of orthogonality 
remains at ${\cal O}(\eps)$ in all of the examples, 
implying backward stability according to
Drko\v{s}ov\'{a} et al.~\cite{drsb1995}.

We anticipate that the iterated Gauss-Seidel algorithms will facilitate the construction of
backward stable iterative solvers based on block Gram-Schmidt algorithms 
including enlarged Krylov subspace methods and $s$-step methods. 
Mixed-precision formulations may also benefit. Randomization and sketching
may be applied to the normal equations in the projector, with
an oblique inner-product $\vect{x}^TB^TB\vect{x}$ and sketching matrix $A$,
leading to even greater computational efficiency. In addition,
these developments could be relevant for Anderson acceleration
of nonlinear fixed-point iterations, which is currently being
applied to optimization algorithms in
deep learning for artificial intelligence. The algorithms may also
be useful for computing eigenvalues with the Krylov-Schur 
algorithm of Stewart~\cite{Stewart86}, where 
${\cal O}(\eps)$ orthogonality is required to obtain an
invariant subspace of maximal dimension.
We plan to explore the parallel strong-scaling performance of the 
IGS-GMRES and hybrid MGS-CGS GMRES for large-scale scientific computation
on exascale supercomputers.

\section*{Acknowledgement}
This research was supported by the Exascale Computing Project (17-SC-20-SC), 
a collaborative effort of the U.S. Department of Energy Office of Science 
and the National Nuclear Security Administration.
The second author was additionally supported by Charles University PRIMUS project no. 
PRIMUS/19/SCI/11 and Charles University Research program no. UNCE/SCI/023. 
In addition, ERC Starting Grant No. 101075632.
The third author was supported by the Academy of Sciences of the Czech
Republic (RVO 67985840) and by the Grant Agency of the Czech Republic,
Grant No. 23-06159S.
Our work was inspired by the contributions 
of A. Bj\"orck, A. Ruhe, C. C. Paige, Z. Strako\v{s}, 
and the legacy of Y. Saad and M. H. Schultz that is MGS-GMRES.
We also wish to thank our friends and colleagues
Julien Langou and Luc Giraud for their thoughtful insights 
over many years. A sincere thank-you to the organizers of the
Copper Mountain conferences on Krylov and multigrid methods for
their rich 38+ year tradition of innovative numerical linear algebra.
We would also like to thank the two anonymous reviewers for their careful 
reading of the manuscript. Their 
thoughtful comments resulted in a much improved presentation.

%

%
%
\newpage

\clearpage
\bibliographystyle{elsarticle-num}
\bibliography{bibGMRES}

\begin{thebibliography}{10}
\expandafter\ifx\csname url\endcsname\relax
  \def\url#1{\texttt{#1}}\fi
\expandafter\ifx\csname urlprefix\endcsname\relax\def\urlprefix{URL }\fi
\expandafter\ifx\csname href\endcsname\relax
  \def\href#1#2{#2} \def\path#1{#1}\fi

\bibitem{Saad86}
Y.~Saad, M.~H. Schultz, {GMRES}: A generalized minimal residual algorithm for
  solving nonsymmetric linear systems, SIAM Journal on scientific and
  statistical computing 7~(3) (1986) 856--869.

\bibitem{2006--simax--paige-rozloznik-strakos}
C.~C. Paige, M.~Rozlo{\v{z}}n{\'{\i}}k, Z.~Strako{\v s}, Modified
  {Gram-Schmidt} ({MGS}), least squares, and backward stability of {MGS-GMRES},
  SIAM Journal on Matrix Analysis and Applications 28~(1) (2006) 264--284.

\bibitem{Paige2002}
C.~C. Paige, Z.~Strako\v{s}, Residual and backward error bounds in minimum
  residual {Krylov} subspace methods, SIAM Journal on Scientific and
  Statistical Computing 23~(6) (2002) 1899--1924.

\bibitem{2020-swirydowicz-nlawa}
K.~\'{S}wirydowicz, J.~Langou, S.~Ananthan, U.~Yang, S.~Thomas, Low
  synchronization {Gram-Schmidt} and generalized minimal residual algorithms,
  Numerical Linear Algebra with Applications 28 (2020) 1--20.

\bibitem{DCGS2Arnoldi}
D.~Bielich, J.~Langou, S.~Thomas, K.~\'{S}wirydowicz, I.~Yamazaki, E.~Boman,
  Low-synch {Gram-Schmidt} with delayed reorthogonalization for {Krylov}
  solvers, Parallel Computing (2021).

\bibitem{Ruhe83}
A.~Ruhe, Numerical aspects of {Gram-Schmidt} orthogonalization of vectors,
  Linear Algebra and its Applications 52 (1983) 591--601.

\bibitem{Bjorck94}
{\AA}.~Bj{\"o}rck, Numerics of {Gram-Schmidt} orthogonalization, Linear Algebra
  and Its Applications 197 (1994) 297--316.

\bibitem{drsb1995}
J.~Drko\v{s}ov\'{a}, A.~Greenbaum, M.~Rozlo\v{z}n\'{i}k, Z.~Strako\v{s},
  Numerical stability of {GMRES}, BIT Numerical Mathematics 35 (1995) 309--330.

\bibitem{Walker88}
H.~F. Walker, Implementation of the {GMRES} method using {Householder}
  transformations, SIAM Journal on Scientific and Statistical Computing 9~(1)
  (1988) 152--163.

\bibitem{gigl:simax:04}
L.~Giraud, S.~Gratton, J.~Langou, A rank-$k$ update procedure for
  reorthogonalizing the orthogonal factor from modified {G}ram-{S}chmidt, SIAM
  J. Matrix Analysis and Applications 25~(4) (2004) 1163--1177.

\bibitem{Bjorck92}
{\AA}.~Bj{\"o}rck, C.~C. Paige, Loss and recapture of orthogonality in the
  modified {Gram-Schmidt} algorithm, SIAM Journal on Matrix Analysis and
  Applications 13 (1992) 176--190.

\bibitem{Lockhart2022}
S.~Lockhart, D.~J. Gardner, C.~S. Woodward, S.~Thomas, L.~N. Olson, Performance
  of low synchronization orthogonalization methods in {Anderson} accelerated
  fixed point solvers, in: Proceedings of the 2022 SIAM Conference on Parallel
  Processing for Scientific Computing (PP), 2022, pp. 49--59.

\bibitem{Bjork_1967}
A.~Bj{\"o}rck, Solving least squares problems by {Gram–Schmidt}
  orthogonalization, BIT 7 (1967) 1--21.

\bibitem{Giraud2005}
L.~Giraud, J.~Langou, M.~Rozlo{\v{z}}n{\'{\i}}k, J.~v.~d. Eshof, Rounding error
  analysis of the classical {Gram-Schmidt} orthogonalization process,
  Numerische Mathematik 101 (2005) 87--100.

\bibitem{Higham1989}
N.~J. Higham, The accuracy of solutions to triangular systems, SIAM J. Numer.
  Anal. 26~(5) (1989) 1252--1265.

\bibitem{HighamBook2002}
N.~J. Higham, Accuracy and Stability of Numerical Algorithms, 2nd Edition,
  SIAM, 2002.

\bibitem{Zhou2021}
Q.~Zou, {GMRES} algorithms over 35 years (2021).
\newblock \href {http://arxiv.org/abs/2110.04017} {\path{arXiv:2110.04017}}.

\bibitem{2020-yamazaki-proceedings-of-siam-pp20}
I.~Yamazaki, S.~Thomas, M.~Hoemmen, E.~G. Boman, K.~{\' S}wirydowicz, J.~J.
  Elliott, Low-synchronization orthogonalization schemes for $s$-step and
  pipelined {Krylov} solvers in {Trilinos}, in: SIAM 2020 Conference on
  Parallel Processing for Scientific Computing, SIAM, 2020, pp. 118--128.

\bibitem{Carson2022}
E.~Carson, K.~Lund, M.~Rozlo{\v{z}}n{\'{\i}}k, S.~Thomas, Block {Gram-Schmidt}
  algorithms and their stability properties, Linear Algebra and its
  Applications 638 (2022) 150--195.

\bibitem{CarsonRoz2021}
E.~Carson, K.~Lund, M.~Rozlo{\v{z}}n{\'{\i}}k, The stability of block variants
  of classical {Gram--Schmidt}, SIAM Journal on Matrix Analysis and
  Applications 42~(3) (2021) 1365--1380.

\bibitem{Stewart86}
G.~W. Stewart, A {K}rylov--{S}chur algorithm for large eigenproblems, SIAM
  Journal on Matrix Analysis and Applications 23~(3) (2001) 601–614.

\bibitem{BienzGroppOlson2019}
A.~Bienz, W.~Gropp, L.~Olson, Node-aware improvements to allreduce, in:
  Proceedings of the 2019 IEEE/ACM Workshop on Exascale MPI (ExaMPI),
  Association for Computing Machinery, 2019, pp. 1--10.

\bibitem{Henrici1962}
P.~Henrici, Bounds for iterates, inverses, spectral variation and fields of
  values of non-normal matrices, Numerische Mathematik 4~(1) (1962) 24--40.

\bibitem{Ipsen1998}
I.~C. Ipsen, A note on the field of values of non-normal matrices, Tech. rep.,
  North Carolina State University. Center for Research in Scientific
  Computation (1998).

\bibitem{Trefethen2005}
L.~Trefethen, M.~Embree, The behavior of nonnormal matrices and operators,
  Spectra and Pseudospectra (2005).

\bibitem{2003_simoncini_SISC}
V.~Simoncini, D.~B. Szyld, Theory of inexact krylov subspace methods and
  applications to scientific computing, SIAM Journal on Scientific Computing
  25~(2) (2003) 454--477.

\bibitem{Green92}
A.~Greenbaum, M.~Rozlo{\v{z}}n{\'{\i}}k, Z.~Strako{\v s}, Numerical behaviour
  of the modified {Gram-Schmidt GMRES} implementation, BIT 37~(3) (1997)
  706--719.

\bibitem{Liesen2004}
J.~Liesen, P.~Tich\'{y}, The worst-case {GMRES} for normal matrices, BIT Numer
  Math 44 (2004) 79--98.

\bibitem{RozloznikST96}
M.~Rozlo{\v{z}}n{\'{\i}}k, Z.~Strako\v{s}, M.~T\accent23uma, On the role of
  orthogonality in the {GMRES} method, in: K.~G. Jeffery, J.~Kr{\'{a}}l,
  M.~Bartosek (Eds.), {SOFSEM} '96: Theory and Practice of Informatics,
  Proceedings, Vol. 1175 of Lecture Notes in Computer Science, Springer, 1996,
  pp. 409--416.

\bibitem{Embree99}
M.~Embree, How descriptive are {GMRES} convergence bounds?, Tech. Rep. Tech.
  Rep. 99/08, Mathematical Institute, University of Oxford, UK (1999).

\bibitem{Mullowney2021}
P.~Mullowney, R.~Li, S.~Thomas, S.~Ananthan, A.~Sharma, A.~Williams, J.~Rood,
  M.~A. Sprague, Preparing an incompressible-flow fluid dynamics code for
  exascale-class wind energy simulations, in: Proceedings of the ACM/IEEE
  Supercomputing 2021 Conference, ACM, 2021, pp. 1--11.

\bibitem{Falgout2004}
R.~D. Falgout, J.~E. Jones, U.~M. Yang, The design and implementation of hypre,
  a library of parallel high performance preconditioners, in: A.~M. Bruaset,
  A.~Tveito (Eds.), Numerical Solution of Partial Differential Equations on
  Parallel Computers, Springer Berlin Heidelberg, Berlin, Heidelberg, 2006, pp.
  267--294.

\bibitem{CHOW2018219}
E.~Chow, H.~Anzt, J.~Scott, J.~Dongarra, Using jacobi iterations and blocking
  for solving sparse triangular systems in incomplete factorization
  preconditioning, Journal of Parallel and Distributed Computing 119 (2018)
  219--230.

\end{thebibliography}

\end{document}